\mag=\magstep1
\documentclass[10pt,a4paper]{amsart}
\usepackage{amssymb,latexsym,amssymb,amscd}

\def\newfont{\bold}
\def\Cal{\mathcal}
\def\frak{\mathfrak}

\def\<<{\langle}
\def\>>{\rangle}

\numberwithin{equation}{section}
\newtheorem{theorem}{Theorem}[section]
\newtheorem{proposition}[theorem]{Proposition}
\newtheorem{corollary}[theorem]{Corollary}
\newtheorem{definition}[theorem]{Definition}

\newtheorem{remark}[theorem]{Remark}
\newtheorem{lemma}[theorem]{Lemma}

\begin{document}

\title{
Theta constants associated to cubic three folds
}

\author{Keiji Matsumoto}
\address{Division of Mathematics, Graduate School of Science, 
Hokkaido University, Japan}
\email{matsu@math.sci.hokudai.ac.jp}

\author{Tomohide Terasoma}
\address{Department of Mathematical Science, University of Tokyo,
Komaba, Meguro, Japan}
\email{terasoma@ms.u-tokyo.ac.jp}

\thanks{The authors would like to express their thanks to Professors  
H. Shiga, M. Yoshida and K. Yoshikawa  for stimulating discussions. 
Professor van Geemen kindly informed us of a conjectural relation between 
our expression of the inverse period map and that of \cite{AF}, see \cite{G}.
}

%    General info
%\subjclass{Primary ; 
% Secondary
% }
\date{July 20, 2000}
%
%\keywords{}

\maketitle

\makeatletter
\renewcommand{\@evenhead}{\tiny \thepage \hfill  K.MATSUMOTO and T.TERASOMA 
\hfill}
\renewcommand{\@oddhead}{\tiny \hfill  THETA FUNCTION OF CUBIC THREE FOLDS
 \hfill \thepage}
\makeatother

\section{Introduction}
\label{sec:introduction}

Every elliptic curve is isomorphic to the double covering 
of the projective line $\bold P^1$ branching at four points. 
The set of the isomorphism classes of elliptic curves 
with marking of branching points can be identified with the configuration 
space $\Cal M_{4pts}$ of ordered four points on $\bold P^1.$ 
By considering their periods, 
one obtains a morphism $per$ from $\Cal M_{4pts}$ to 
the set $\Cal M_{ct}$ of the isomorphism classes of complex torus 
with marking of $2$-torsion points.  
The set $\Cal M_{ct}$ can be identified with the quotient 
space ${\frak H}/\Gamma(2)$ 
of the upper half plane $\frak H$ by 
the principal congruence subgroup $\Gamma(2)$ of level $2$ 
in $SL(2,\bold Z).$ 
By the classical theory of elliptic functions,
the map $per : \Cal M_{4pts} \to \Cal M_{ct}$ is an isomorphism.
The inverse of $per$ can be described in terms of 
celebrated Jacobi's theta constants. 

Many mathematicians have tried to find moduli spaces of suitable algebraic varieties 
which can be uniformized by some symmetric space.  
E. Picard was the first person who found such moduli space 
which is two dimensional.  
He studied a family of Picard curves, which are cyclic triple coverings 
of $\bold P^1$ 
branching at five points. 
The periods of a curve determine an element of
the $2$-dimensional complex ball $\bold B_2$ embedded 
in the Siegel upper-half space $\frak H_3$ of degree $3$.
This correspondence gives a uniformization of 
the moduli space of the Picard curves by $\bold B_2$.
The inverse $\Theta$ of the period map was expressed 
in terms of special values of theta functions for 
the Jacobians of Picard curves.
Shiga then found the representation of the map $\Theta$  in 
terms of theta constants, which are automorphic forms on 
$\bold B_2$ with respect to the monodromy group. 
Inspired by Picard's results, 
the moduli spaces of the cyclic coverings of $\bold P^1$ branching at 
$(n\!+\!3)$-points uniformized by the $n$-dimensional complex balls 
were classified by Terada (\cite{T}), Deligne and Mostow (\cite{DM}).
One specific  three dimensional moduli spaces listed in \cite{DM} 
was studied in \cite{Ma} similarly to Shiga: 
the inverse of the period map for a family of cyclic triple coverings of 
$\bold P^1$ branching at six points was expressed 
in terms of theta constants.

Recently, Allcock, Carlson and Toledo showed that 
the moduli space of marked cubic surfaces can be uniformized by 
the $4$-dimensional complex ball $\bold B_4$, 
though the period map for any family of cubic surfaces is constant.
Let $Y$ be the cyclic triple covering of 
$\bold P^3$ branching along a cubic surface $X$. 
The intermediate Jacobian $J(Y)$ of $Y$ is a $5$-dimensional abelian variety. 
By considering the normalized period matrix of $J(Y),$ 
they obtain a point $\tau$ in the Siegel upper half space 
$\frak H_5$ of degree $5$.
Since the abelian variety $J(Y)$ admits an action of $\mu_3$, 
$\tau$ belongs to the subdomain 
$\bold B_4=\{\tau\in \frak H_5\mid (H\tau)^2+H\tau+I=0\}$ in $\frak H_5$, 
where $H={\rm diag}(1,1,1,1,-1)$ and 
$\mu_m$ is the group of $m$-th roots of unity.
As a consequence, they get a multivalued holomorphic
map $\varphi$ from 
the moduli space $\Cal M_{cs}$ of cubic surfaces with marking of the $27$ lines 
to the 4 dimensional complex ball $\bold B_4$.  
Its image is an analytic Zariski open set of $\bold B_4$.  
In this manner, they get a period map.
Moreover, the monodromy group of $\varphi$ is an arithmetic subgroup $\Gamma$ 
of the unitary group $U(4,1)$ and the induced holomorphic map 
$per : \Cal M_{cs} \to \Gamma\backslash \bold B_4$ is a birational
morphism.  They define an action of $PO(5, \bold F_3)\simeq W(E_6) $ on
$\Gamma\backslash \bold B_4$ which is compatible with
the classical action of the Weyl group $W(E_6)$ on $\Cal M_{cs}$
through the period map $per$. 
	
In this paper, 
we study the action of $W(E_6)$ on the 
$(1-\rho)$-torsion subgroup $J(Y)_{1-\rho}(\simeq \bold F_3^5)$ of $J(Y),$ 
where $\rho$ is an automorphism of $Y$ of order $3$. 
The intersection form on 
$H^{prim}_2(X,\bold Z)$ defines a $\bold F_3$-valued quadratic form $q$  
invariant under the action of $W(E_6)$. 
To each element $v$ in $J(Y)_{1-\rho}$, we assign a theta constant 
$\Theta_v(\tau)$ on $\bold B_4.$ 
It is easy to see that $\Theta_v(\tau)\equiv 0$ for $v\notin 
S=\{v\in J(Y)_{1-\rho} \mid q(v)= 0\}.$ 
The cubes of the non vanishing eighty $(=\# S)$ theta constants 
for $v\in S$ give a $W(E_6)$-equivariant projective 
embedding $\Theta$ of $\Gamma\backslash \bold B_4$.
On the other hand, the moduli space of smooth cubic surfaces with 
marking of the $27$ lines is isomorphic to the 
moduli space $\Cal M_{6pts}$ of ordered 
$6$ points on $\bold P^2$ in general position. 
After Coble, there exists a $W(E_6)$-equivariant 
projective embedding $Z$ of $\Cal M_{6pts}$ using $80$ polynomials labeled by $S$.  
The main theorem (Theorem \ref{main theorem}) of this paper
asserts that the projective embeddings $\Theta$ 
and $Z$ coincide via the period map $\varphi$.
As a consequence, we obtain the inverse of the period map $\varphi.$

Here we explain the basic tool for our study. For a line $L$ in $X$, 
we associate a $\mu_3$-covering $C$ of $L$ branching at twelve points. 
This curve $C$ admits an involution $\sigma$ on
$C$ commuting with $\rho$ and fixing exactly two points among
the twelve,
which are denote by
$\{p_0 = \sigma (p_0), p_\infty= \sigma (p_\infty ),
p_1,\dots ,p_5, \sigma (p_1), \dots ,\sigma (p_5)\}$.
The key fact in this paper is that the Prym variety 
$Prym(C,\sigma)$ is isomorphic to $J(Y)$.  
As a consequence, the periods of $J(Y)$ are equal to
those of $Prym(C, \sigma)$. 
Note that the brancing index of the $\mu_6$-covering 
$C\to C/<\sigma,\rho>\simeq \bold P^1$ 
appears in the list of Mostow's paper \cite{Mo}. 
A proof of the surjectivity of the period map after Allcock, Carlson and 
Toledo  reduces to Mostow's results.
Via the isomorphism $Prym(C,\sigma)\simeq J(Y)$, 
the images of $p_1, \dots ,p_5$ under the Abel-Jacobi map 
$C \to Prym(C, \sigma )$
form an orthonormal basis in 
$J(Y)_{1-\rho}$.
Conversely any orthonormal basis of $J(Y)_{1-\rho}$ can be obtained from
one of the $27$ lines.  Thus we have a one to one correspondence
between the $27$ lines and orthonormal bases of $J(Y)_{1-\rho}$.
This correspondence yields a dictionary between  
the geometry of the $27$ lines and the geometry of $(\bold F_3^5,q)$.
For example, two lines $L_1$ and $L_2$ intersect if and only
if the corresponding orthonormal bases contain a common element
up to signature.  From this fact, we get a one to one correspondence
between the $45$ tritangents and the elements of length $1$, modulo signature.
Moreover we get a one to one correspondence between the $36$ double sixes and
the elements of length $2$, modulo signature. 
 As an application of this dictionary, we describe
relations of degree $3$ and $9$ for theta constants in terms of
the geometry of $(\bold F_3^5,q)$.

We would like to explain the contents of this paper.
In Section 2, we study the line geometry on the cubic threefold $Y$
obtained by the $\mu_3$-covering of $\bold P^3$ branching along
a cubic surface $X$.  There we introduce a curve $C$, 
an involution $\sigma$ and an automorphism $\rho$ of order three.
The main result in this section is Corollary \ref{cor:prym}:
the Prym variety of $C$ for the involution $\sigma$ is
isomorphic to the intermediate Jacobian $J(Y)$ of $Y$.

In Section 3, we introduce a symplectic basis 
$A_1,\cdots,A_5,B_1,\dots,B_5$ of
the intermediate Jacobian $J(Y)$ of $Y$, which will be used to
describe theta functions in Section 4.
Thanks to the explicit description of the basis of topological 
cycles, we check that the images of the branch points under 
the Abel-Jacobi map form an orthonormal basis of the
$(1-\rho)$-torsion subgroup $J(Y)_{1-\rho}$ of $J(Y)$.
There we describe the correspondence between the $45$ tritangents and 
the elements of length 1 in $\bold F_3^5$ and an explicit isomorphism
between $PO(5,\bold F_3)$ and $Aut (\Gamma_{std})$,
where $\Gamma_{std}$ is the dual graph of a cubic surface.

We introduce the theta function associated to the Prym variety
equipped with a symplectic basis.
We study the pull-back of theta functions  
by morphisms $C\to Prym (C,\sigma)$.

In the last section, we state the main theorem.
We define
a projective embedding $Z:\Cal M_{6pts} \to \bold P^{79}$ of  
$\Cal M_{6pts}$ by $80$ polynomials 
labeled by $S. $
The advantage of using these polynomials is that they behave well
under the action of $W(E_6)$.  Using this description, we prove
that the inverse of the period map can be described by theta 
constants defined in Section 4.

Before closing this introduction, we give some remarks on 
degenerations of cubic surfaces.
A nodal cubic surface is obtained from the blowing up
of six points on a conic $C$.
The intermediate Jacobian of the $\mu_3$-covering of the
nodal cubic surface is isomorphic to the Jacobian of
a $\mu_3$-covering $\tilde C$ of the conic $C$ 
branching at the six points.
Actually, the genus of this curve is $4$ and the image 
of the canonical map to $\bold P^3$
is a complete intersection of hypersurfaces 
$H_2$ and $H_3$ of degree 2 and 3 in $\bold P^3$.
Let $\tilde P$ be the blowing up 
of $\bold P^3$ along the embedded curve $\tilde C$. Then we can
contract the strict transform of $H_2$.  
The resulting threefold $Y$ is nothing
but the $\mu_3$-covering of $\bold P^3$ branching along the
nodal cubic surface.
In the closure $\bar \Cal M_{cs}$ of the image of $Z:\Cal M_{cs}\to\bold P^{79}$,
the component of the boundary $\bar \Cal M_{cs} - \Cal M_{cs}$
corresponds to the mirrors of the 36 reflections in $W(E_6)$.
The detail will appear elsewhere.

{\bf Notations.}

\begin{itemize}
\item[$\mu_m$]:  The group of $m$-th roots of unity in $\bold C^{\times}$. 
\item[$\Im (\tau)$] :  The imaginary part of a complex matrix $\tau$.
\item[$M_0$] :  The row vector consisting of the diagonal entries
of an $n\times n$ matrix \\ \hspace{0.1cm} $M$. 
\item[$\bold e (x)$] $:= \exp ( 2\pi \sqrt{-1}x)$. 
\item[$\bold 1$] $:= (1,\dots ,1)$. 
\end{itemize}

\section{Line geometry of cubic threefold with $\mu_3$-actions}
\label{sec:linegeometry}

\subsection{A normal form for cubic surfaces}
\label{subsec:normal form}
In this section, we introduce a certain normal form
for cubic surfaces.
First we recall several well known facts about smooth cubic surfaces.
Let $X$ be a cubic surface, i.e. a smooth surface of degree
3 in $\bold P^3$. There are exactly 27 lines on
the surface $X$ and we can choose
mutually disjoint 6 lines from them,
which are written as $E_1, \dots ,E_6$. Since the self-intersection
numbers
of these lines are $-1$, by contracting $E_1, \dots , E_6$ to
points $P_1, \dots ,P_6$, we obtain
the two dimensional projective space $\bold P^2$. 
The set of points $P_1, \dots , P_6$ on $\bold P^2$ 
are generic in the following sense:
\begin{enumerate}
\item No three of $P_1, \dots ,P_6$ are collinear.
i.e. there exist no lines
passing through three points among the six.
\item The six points $P_1, \dots ,P_6$ are not coconic, 
i.e. there exist no conics passing through
the six.
\end{enumerate}
For $i=1, \dots , 6$, there exists a unique conic $\bar C_i$ 
(resp. a line $\bar L_{ij}$)
in $\bold P^2$ containing
the $P_j$'s $(j \neq i)$ (resp. $P_i$ and $P_j$ ($1 \leq i < j \leq 6$)).  
The proper transforms $C_i$ and 
$L_{ij}$ of $\bar C_i$ and $\bar L_{ij}$
are lines in $X$.
Then $E_1, \dots ,E_6$, $C_1, \dots ,C_6$ and $L_{ij}$
$(1 \leq i < j \leq 6)$ are the 27 distinct lines in the cubic surface
$X$. 

A notion of a marked cubic surface is defined as follows.
We define the standard dual graph $\Gamma_{std}$ as follows.
(1) The set of the vertices consists of $e_i, c_i$ $(i=1, \dots , 6)$
and $l_{ij}$ $(1\leq i<j \leq 6)$.
(2) $e_i$ is adjacent to $c_j$ if and only if $i \neq j$.
(3) $e_i$ (resp. $c_i$) is adjacent to $l_{jk}$ if and only if $i \in \{j,k\}$.
(4) $l_{ij}$ is adjacent to $l_{kl}$ if and only if $\{ i,j\} \cap
\{ k,l\}= \emptyset$.
(5) $e_i$ and $e_j$ (resp. $c_i$ and $c_j$) ($i \neq j$) 
does not adjacent to each other.
Then the dual graph $\Gamma (X)$ of the 27 lines in $X$ 
is isomorphic to $\Gamma_{std}$. A marking $\Psi_{cs}$ 
of $X$ is defined as an isomorphism 
$\Psi_{cs} : \Gamma (X) \to \Gamma_{std}$ of graphs.

Since $P_6$ is not contained in the conic
$\bar C_6$ in $\bold P^2$, there are two tangent lines
$T_0$ and $T_{\infty}$
of $\bar C_6$ passing through $P_6$.  
The tangent points are denoted by $Q_0$ and
$Q_{\infty}$.
If $\bar L_{i6}$ tangents to $\bar C_6$ in $\bold P^2$, 
then $L_{i6}\cap C_6\cap E_i$ 
consists of one point $P$.
This point $P$ is called an Eckardt point of $X$. 
If a cubic surface is sufficiently
generic, then there exist no Eckardt points.
Until the end of this section, we assume that $X$ has no Eckardt points.
We choose a coordinate $(x_0:x_1:x_2)$ of $\bold P^2$
such that 
\begin{enumerate}
\item
the conic $\bar C_6$ is expressed as $x_0x_1=x_2^2$, and
\item
the equations of the two tangent lines $T_0$ and $T_{\infty}$ are 
$x_0=0$ and $x_1=0$, respectively.
\end{enumerate}
Then $P_i$ $(i=1, \dots , 5)$ are given by
$\displaystyle (\frac{1}{a_i}: a_i :1)$ ($a_i \neq 0$)
and $P_6$ by $(0:0:1)$.
Now we define a polynomial $h(x)$ and $s_i$ $(i=1, \dots ,5)$ 
by
\begin{align*}
h(x) = &  \prod_{i=1}^5(x - a_i) \\
= & x^5+ s_1x^4+s_2x^3+ s_3x^2+s_4x+s_5. 
\end{align*}
The vector space of homogeneous polynomials of degree three  
on $x_0, x_1, x_2$ vanishing at
6 points $P_1, \dots ,P_6$ is four dimensional.
A basis of this space is given by
\begin{align*}
u_0 &=(x_0x_1-x_2^2)x_0, \\
u_1 &=(x_0x_1-x_2^2)x_1, \\
u_2 &=x_1^3+s_1x_1^2x_2+ s_2x_1x_2^2 + s_3x_0x_1x_2 
+ s_4 x_0x_2^2 +s_5 x_0^2x_2, \\
u_3 &=x_1^2x_2+s_1x_1x_2^2+ s_2x_0x_1x_2 + s_3 x_0x_2^2 
+ s_4 x_0^2x_2 +s_5 x_0^3.
\end{align*}
By eliminating $x_0,x_1,x_2$, we get the cubic relation 
$F(u_0, u_1, u_2, u_3)=0$, where
\begin{align*}
F(u_0, u_1, u_2, u_3) =& u_0u_2^2+(s_2u_0u_1+s_4u_0^2 -u_1^2)u_2 \\
& -(u_1u_3^2+(s_3u_0u_1+ s_1u_1^2-s_5u_0^2)u_3) \\
& -( s_2u_1^3 + s_4u_1^2u_0-s_1s_5u_0^2u_1 -s_3s_5u_0^3). 
\end{align*}
This is the defining equation of the cubic surface $X$.  
Under these coordinates $(u_0:\cdots :u_3)$, the proper
transform $C_6$ of the conic $\bar C_6$ is given by $u_0 = u_1= 0$. 
Any point of $\bar C_6$ different from $Q_0$ and $Q_{\infty}$ is
expressed as $\displaystyle (\frac{1}{x}:x:1)$,  $(x \neq 0)$
and the corresponding point in the 
proper transform $C_6$ is given by $u_0 = u_1 =0$, 
$u_2 = \displaystyle \frac{1}{x^2}h(x)$ and
$u_3 = \displaystyle \frac{1}{x^3}h(x)$. The projective coordinates of
this point are given as
$$
(u_0 : u_1 : u_2 : u_3) = (0:0:x:1).
$$
\subsection{Triple covering of $\bold P^3$ branching along the cubic surface}
\label{subsec:triple covering}
Let $F=F(u_0, u_1, u_2, u_3)$ be the defining equation of 
the cubic surface $X$ as in \S \ref{subsec:normal form}.
The cubic threefold $Y$ in $\bold P^4$
defined by $u_4^3=F(u_0,u_1,u_2,u_3)$ 
is a cyclic covering of $\bold P^3$ branching along 
$X\simeq \{(u_0 :\cdots : u_4)\in Y\ \mid u_4=0\}$.
We define an action $\rho$ of $\zeta \in \mu_3$ on $Y$ by
$$
\rho : (u_0:\cdots :u_3:u_4) \mapsto (u_0: \cdots : u_3 :\zeta u_4).
$$ 
\begin{definition}
\label{def:fano var}
The subvariety of the Grassmann variety $Gr (5,2)$
of lines in $\bold P^4$ consisting of
lines in $Y$ is called the Fano variety $F(Y)$ of $Y$. (see \cite{CG}.)
For a line $l$ in $Y$, the subvariety $Inc(l)$ of $F(Y)$ consisting of lines
which intersect with $l$ is called the incidental subvariety of $l$.
\end{definition}
\begin{remark}
\label{rem:CG incidental}
It is known (c.f.\cite{CG}) that the Fano variety $F(Y)$ is
a smooth surface and that $Inc(l)$ is a divisor of $F(Y)$.
\end{remark}
Since $X$ is the branch locus of the covering $Y \to \bold P^3$,
a line $L$ in $X$ is a line in $Y$.   
In the rest of this section, we investigate  
the incidental subvariety $Inc(C_6)$ of the line $C_6$ in $Y$
and its normalization $C$.

Let $P$ be a point of $C_6$ in $Y$ different from 
$Q_{\infty},Q_0$.
Then its coordinates are $(u_0:u_1:u_2:u_3:u_4)=(0:0:x:1:0)$ 
($x \neq 0$).
The line $L$ joining two points $(0:0:x:1:0)$ and $(a:b:c:0:e)$ 
can be parameterized as $(u_0: \dots :u_4) = (at, bt, ct+x, 1, et)$, 
$(t \in \bold C \cup \{ \infty \})$.
We consider the condition for $a,b,c,e$ so that
the line $L$ is contained in $Y$, i.e.
$F(at, bt, ct+x, 1)=(et)^3$ holds for all $t$. 
By the straight forward calculation, we have
\begin{align*}
ax^2= & b, \\
2acx+(s_2ab+s_4a^2-b^2)x= & s_3ab+s_1b^2-s_5a^2, \\
ac^2+(s_2ab+s_4a^2-b^2)c = &  (s_2b^3 +s_4b^2a-s_1s_5a^2b-s_3s_5a^3)+e^3. 
\end{align*}
By putting $a=1$, we get $b=x^2$ and
\begin{align*}
c= & \frac{1}{2x}(x^5+s_1x^4-s_2x^3+s_3x^2-s_4x-s_5), \\
4x^2e^3 = & h(x)h(-x). 
\end{align*}
By setting $y=4^{1/3}xe$, we get a family of lines
contained in $Y$ parameterized by the curve $C^0$:
\begin{equation}
\label{eq:curve}
y^3 = x h(x)h(-x)=x\prod_{i=1}^5(a_i^2 -x^2)
\quad (x \neq 0, \infty).
\end{equation}
By the map $C \ni (y,x) \to (\displaystyle\frac{1}{x}:x:1) \in C_6$,
$C$ is regarded as a $\mu_3$-covering of $C_6$.
By direct computation, we can readily show that
there are no lines passing through $Q_0$
(resp. $Q_{\infty}$) other than $C_6$.  Therefore 
$Inc(C_6)=C^0 \cup \{ [C_6]\}$,
where $[C_6]$ is the point of $Gr(5,2)$ corresponding to the line
$C_6$.  
Since $Inc(C_6)$ has two tangents at $[C_6]$ in $Gr(5,2)$,
it is a curve only with a node.  
Thus we get a family of lines in $Y$ parameterized by $C$.

Note that the action of $\mu_3$ on the Fano variety is compatible with
the action of $\zeta \in \mu_3$ on the curve $C$ defined by 
$(y \to \zeta y, x \to x) \in Aut(C)$.  The action of 
$\displaystyle \omega=\frac{-1+\sqrt {-3} }{2}$
is denoted by $\rho$.

\subsection{Cylinder map}
\label{subsec:cylinder map}
  In this section, we introduce and investigate the cylinder map 
induced by the family of lines defined in \S\ref{subsec:triple covering}.
Let $\Cal U$ be the universal family defined by
$\Cal U = \{ (p, x) \mid p \in C, x \in l_p\}$, where $l_p$ is
the line corresponding to the point $p \in C$, and let
$pr_1 :\Cal U \to C$ and $pr_2 :\Cal U \to Y$
be the natural projections. Then $pr_1$ induces an isomorphism:
$$
pr_1^*:H^1(C, \bold Z(-1)) \overset{\simeq}\longrightarrow 
H^1(\Cal U, \bold Z(-1)).
$$
The Gysin map $pr_{2*}$ is the Poincare dual of the natural homomorphism
$pr_2^*:H^3(Y, \bold Z(3)) \to H^3(\Cal U, \bold Z(3))$, i.e.
$pr_{2*}b$ ($b \in H^1(\Cal U, \bold Z(-1))$) is an element so that
\begin{equation}
\label{eq:compat symplectic}
a\cup pr_{2*}b = pr_2^* a \cup b
\end{equation} 
holds for any $a \in H^3(Y,\bold Z)$. 
We define the cylinder map:
$$
c:H^1(C, \bold Z(-1)) \overset{pr_1^*}\longrightarrow 
H^1(\Cal U, \bold Z(-1)) 
\overset{pr_{2*}}\longrightarrow H^3(Y, \bold Z).
$$
This is a homomorphism of Hodge structures.

We define an involution $\sigma$ on $C$ by 
$y \to -y$ and $x \to -x$.
\begin{definition}
\label{def:CG invol}
An involution $\sigma '$ on $C$ is called the Clemens-Griffiths involution if
$l_p$ and $l_{\sigma ' (p)}$ are contained in a plane in $\bold P^4$ for
generic $p \in C$.
\end{definition}
\begin{proposition}
\label{prop:GC-involution}
The involution $\sigma$ coincides with the Clemens-Griffiths involution.
Moreover the action of $\sigma$ commutes with the action of $\rho$.
As a consequence the group $\mu_6$ of 6-th roots of unity acts
on the curve $C$.
\end{proposition}
\begin{proof}
By the definition of $\sigma$, if $l_p$ is a line connecting $(0:0:x:1:0)$
and $(1:b:c:0:e)$, then $l_{\sigma (p)}$ is a line connecting
$(0:0:-x:1:0)$ and $(1:b:c':0:e)$. Therefore both $l_p$ and $l_{\sigma (p)}$
are contained in the plane $bu_0 = u_1, eu_0 = u_4$.
The commutativity of the action of $\rho$ and $\sigma$ is a direct consequence
of the definition of $\sigma$ and $\rho$.
\end{proof}
\begin{corollary}
Let $H^1(C, \bold Z(-1))^+ =\{ v \in H^1(C, \bold Z(-1)) 
\mid \sigma^*(v) = v \}$,
and $H^1(C, \bold Z(-1))^- = H^1(C, \bold Z(-1))/H^1(C, \bold Z(-1))^+$.
Then the cylinder map $c$ factors through $H^1(C, \bold Z(-1))^-$.
\end{corollary}
\begin{proof}
By Proposition \ref{prop:GC-involution},
$H^1(C, \bold Z(-1))^-$ is the maximal quotient
on which the Clemens-Griffiths involution acts as the $(-1)$-multiplication.
On the other hand,
the Clemens-Griffiths involution acts on 
$H^3(Y, \bold Z)$ as the $(-1)$-multiplication 
and we obtain the corollary.
\end{proof}
The induced map $H^1(C, \bold Z(-1))^- \to H^3(Y, \bold Z)$ is denoted
by $\phi$.
The algebra over $\bold Z$ generated by $\rho$ with the relation
$1 + \rho + \rho^2 =0$ is denoted by $\bold Z[\rho ]$. Then
$H^1(C, \bold Z(-1))^-$ and $H^3(Y, \bold Z)$ 
are modules over $\bold Z[\rho ]$.
\begin{theorem}
\label{theorem:prym}
The morphism $\phi$ is an isomorphism as 
$\bold Z[\rho ]$ modules.
\end{theorem}
\begin{proof}
Since the action of $\rho$ is compatible, 
$\phi$ is a homomorphism as $\bold Z[\rho ]$ modules.
We will prove that this
is actually an isomorphism.
Let $D$ be the quotient of $C$ by the involution $\sigma$. Then 
by the Hurwitz theorem, the genus $g(C)$ and $g(D)$ of $C$ and $ D$ are 
$10$ and $5$, respectively. Since the rank of 
$H^3(Y, \bold Z)$ is $10$, it is enough to prove the surjectivity.

We use the same notations $C_6$, $Inc(C_6) \subset Gr(5,2)$ as in
the last paragraph.  Let $g:\Cal  Y \to \Delta$ be
a small deformation of the cubic threefold $\Cal Y_0 =Y$ 
over $\Delta =\{ t \in \bold C \mid \mid t \mid < \epsilon \}$ such
that the generic fibers $\Cal Y_t$ at $t\in \Delta^*=\Delta -\{ 0\}$ are
sufficiently generic.  We can extend a line $C_6$
to a family of lines $\Cal C \to \Delta$ contained in $\Cal Y$.
The family of incidental subvarieties $f:\Cal D \to \Delta$ 
is a family of curves
on $\Delta$. By \cite{CG}, the generic fiber $\Cal D_t$ of $\Cal D$ at 
$t \in \Delta^*$ is a smooth curve of genus 11. Moreover the family of
the Clemens-Griffiths involutions on the fibers comes to be an involution
$\sigma$ of $\Cal D$ preserving each fiber.
Thus we have the relative cylinder map
$$
\bold R^1 f_* \bold Z(-1) \to
\bold R^3 g_* \bold Z,
$$
and it factors through the maximal quotient 
$(\bold R^1 f_* \bold Z(-1))^-$ of
$\bold R^1 f_* \bold Z(-1)$ on which the involution $\sigma$ acts
as the $(-1)$-multiplication:
$$
(\bold R^1 f_* \bold Z(-1))^{-} \to
\bold R^3 g_* \bold Z.
$$
Since $g$ is a proper smooth morphism of cubic threefolds, 
$\bold R^3g_*\bold Z$ is a smooth $\bold Z$-sheaf.
By considering
the rank of fibers and the specialization map, we get the 
smoothness of the sheaf $(\bold R^1 f_* \bold Z(-1))^{-}$.
On the other hand, by Theorem 11.19 of \cite{CG}, for $t \in \Delta^*$
the homomorphism
$$
(\bold R^1 f_* \bold Z(-1))_t^{-} \to
\bold R^3 g_* \bold Z_t
$$
is surjective. Therefore
$$
(\bold R^1 f_* \bold Z(-1))_0^{-} \to
\bold R^3 g_* \bold Z_0
$$
is surjective.  Since $Inc(C_6)$ is a curve with one node,
whose normalization $C$ is a curve of genus 10, we have
$$
(\bold R^1 f_* \bold Z(-1))_0^{-} \simeq H^1(C, \bold Z)^{-},
$$
and we get the theorem.
\end{proof}
\begin{remark}
\label{rem:pol J and C}
By the compatibility of cup products (\ref{eq:compat symplectic}),
the polarization of $J(Y)$ given by the cup product
is equal to the half of the restriction of
the cup product on $H^1(C,\bold Z)$ to $H^1(C, \bold Z)^-$.
\end{remark}
Let $C$, $\sigma : C \to C$ and $D = C/<\sigma >$ be as in 
the proof of Theorem \ref{theorem:prym}.
Since the double covering $C \to D$ branches at two points, the
Prym variety $Prym (C,\sigma) = Ker (J(C) \to J(D))$ 
is a principally polarized abelian variety (c.f.\cite{Mu}).  
This is isomorphic to the image of the morphism
\begin{equation}
\label{eqn:1-sigma}
(1-\sigma) : J(C) \to J(C): v \to v -\sigma (v).
\end{equation}
Via this isomorphism, $Prym(C, \sigma)$ is regarded as the
maximal quotient of $J(C)$ on which the involution $\sigma$ acts
as $(-1)$-multiplication.
We have the following corollary to Theorem \ref{theorem:prym}.
\begin{corollary}
\label{cor:prym}
Let $cor_C:J(C) \to J(Y)$ be the homomorphism induced
by the correspondence $C \leftarrow \Cal U \to Y$,
where $J(Y)$ is the intermediate Jacobian of $Y$.
The map $\bar c:Prym(C, \sigma) \simeq J(Y)$ induced by $-cor_C$
is an isomorphism.
\end{corollary}
By the commutativity of $\sigma$ and $\rho$, $\rho$ induces 
an automorphism of $Prym(C, \sigma)$. The induced automorphism 
is also denoted by $\rho$.
It is easy to see that $\bar c$ is
compatible with the action of $\mu_3=< \rho >$.

\subsection{The Abel-Jacobi map and the level map $\Lambda$}
\label{subsec:abel-Jacobi and orthogonal}
We consider the following commutative diagram:

\setlength{\unitlength}{0.75mm}
\begin{picture}(160,60)(-17,-5)
\put(0,45){$H_2^{prim}(X, \bold Z)$}
\put(50,45){$H_2(X, \bold Z)$}
\put(100,45){$H_2(\bold P^3, \bold Z)$}
\put(8,30){$\lambda$}
\put(8,8){$a$}
\put(50,30){$CH_1(X)$}
\put(100,30){$CH_1(\bold P^3)$}
\put(0,15){$A_1(Y)$}
\put(50,15){$CH_1(Y)$}
\put(100,15){$CH_1(\bold P^4)$}
\put(0,0){$J(Y)$}
\put(53,37){$\simeq$}
\put(103,37){$\simeq$}
\put(30,46){\vector(1,0){17}}
\put(73,46){\vector(1,0){22}}
\put(73,31){\vector(1,0){22}}
\put(15,16){\vector(1,0){32}}
\put(73,16){\vector(1,0){22}}
\put(5,40){\vector(0,-1){20}}
\put(60,35){\vector(0,1){8}}
\put(60,28){\vector(0,-1){8}}
\put(110,35){\vector(0,1){8}}
\put(110,28){\vector(0,-1){8}}
\put(5,13){\vector(0,-1){8}}
\end{picture}
where $CH_1(X)$ is the Chow group of $X$ of dimension 1, $A_1(Y)$
is the subgroup of $CH_1(Y)$ consisting of algebraic cycles
algebraically equivalent to zero, and $H_2^{prim}(X, \bold Z)$
is the kernel of the natural homomorphism 
$H_2(X, \bold Z) \to H_2(\bold P^3, \bold Z)$. 
Let $l$ and $l'$ be lines in $X$ disjoint.
Then the images of lines $l$ and $l'$ in $Y$ are
algebraically equivalent. Thus the image of $[l]-[l']$
in $CH_1(Y)$ is contained in $A_1(Y)$.
Since $H_2^{prim}(X, \bold Z)$ is generated by elements $[l]-[l']$,
with $l \cap l' = \emptyset$, the image of $H_2^{prim}(X,\bold Z)$ is
contained in $A_1(Y)$. The induced map 
$H_2^{prim}(X, \bold Z) \to A_1(Y)$
is denoted by $\lambda$. By composing $\lambda$ and the Abel-Jacobi map
$a:A_1(Y) \to J(Y)$, we get a
homomorphism $\Lambda =a \circ \lambda :H_2^{prim}(X, \bold Z) \to J(Y)$,
which is called the level map for $X$.

Let $X$ be a cubic surface with a marking $\Psi_{cs}$. 
We define an element $v_i$ in $J(Y)$ by 
$$
v_i=\Lambda([E_i]-[L_{i6}]).
$$ 
The point of $C$ defined by $x=a_i, y=0$ (resp. $x=0, \infty $) 
in the equation (\ref{eq:curve})
is denoted as $p_i$ (resp. $p_0, p_\infty$).
The point $p_i$ (resp. $\sigma (p_i)$) corresponds to the line
$E_i$ (resp. $L_{i6}$).  
We define two morphisms $jac:C\to J(C)$ and $\jmath:C\to Prym(C,\sigma)$:
\begin{align*}
C \ni p &\mapsto [p]-[p_0]\in J(C),\\
C \ni p &\mapsto [p]-[\sigma(p)]\in Prym(C,\sigma),
\end{align*}
where $[p]$ is the divisor class of $p$. 
Then we have the following commutative diagram:
\begin{equation}
\label{jac and jmath}
\matrix
C  &  &\overset \jmath\longrightarrow & &Prym(C)\\
   &  &                & &       \\
jac \downarrow\qquad &  &\footnotesize{(1-\sigma)}\nearrow\qquad\quad& 
&\quad \downarrow \bar c\\
   &  &                & &       \\
J(C)& &\underset{-cor_C}\longrightarrow & &J(Y).  \\
\endmatrix
\end{equation}
\begin{proposition}
\label{prop:cor vi and ai}
Under the composite map $\bar c\circ \jmath$
$$C\overset\jmath \longrightarrow Prym(C,\sigma) \overset {\bar c}\simeq J(Y),
$$
the image of $p_i \in C$ is $v_i$.
\end{proposition}
\begin{proof}
By the definition of $jac$ and $\bar c$, we have  
$$jac(p)=[p]-[p_0], \quad 
-cor([p]-[p_0])=-\Lambda([E_i]-[C_6]).$$
Since the involution $\sigma$ induces the $(-1)$-multiplication on 
$Prym(C,\sigma),$ we have 
$-\Lambda([E_i]-[C_6])=\Lambda([E_i]-[L_{i6}]).$
The diagram (\ref{jac and jmath}) yields the proposition.
\end{proof}

\section{Finite geometry for  $(1-\rho)$-torsion subgroup of the Prym variety}
\label{sec:finite geometry}

\subsection{Symplectic basis for the curve $C$ and $Prym (C,\sigma)$}
\label{subsec:symplectic basis}
In this section we introduce a symplectic basis of $J(C)$ and 
$Prym( C, \sigma)$
compatible with the action of $\rho$. 
We assume that the cubic surface $X$ has no Eckardt points in this
section.
To specify topological cycles of 
$C$, we assume that $a_1, \dots ,a_5 \in \bold R$ 
and $0< a_1 < \cdots < a_5$.
The curve $C$ can be expressed as $\mu_6$ covering of
$\bold P^1\simeq C/<\rho , \sigma >$:
put $\xi = x^2$ in (\ref{eq:curve}), then the curve $C$ is given by
\begin{equation}
\label{eq:degree 6 cov eq}
y^6=\xi \prod_{i=1}^5(\xi -a_i^2)^2.
\end{equation}
Note that 
$x=\prod_{i=1}^5(a_i^2 - \xi)^{-1}y^{3}$.  The actions of $\rho$ and $\sigma$
are given by
\begin{align*}
\rho (y) = \omega y, \quad \rho (\xi )= \xi, \\
\sigma (y) = - y, \quad \sigma (\xi )= \xi. 
\end{align*}
By gluing 6 copies of the $\xi$-planes cut along the slit
as in Figure 1, we get the curve $C$.
Here the projection of the point $p_i$ 
to $\xi$ plane is denoted by $\xi_i$ for short. 
We define topological cycles
$\beta_1, \dots , \beta_5$ 
as in the Figure 1. In Figure 1 the numbers written along paths are
the label of sheet through which the paths are passing.
Here the sheets are labeled with $\bold Z / 6\bold Z$.
The action of $\sigma$ (resp. $\rho$) sends the $i$-th sheet to the
$(i+3)$-th sheet (resp. $(i+2)$-th sheet).

\setlength{\unitlength}{0.75mm}
\begin{figure}[hbt]
\begin{picture}(160,50)(-5,10)
%cut lines
\thicklines
\put(-5,15){\line(1,0){150}}
\put(-5,55){\line(1,0){150}}
\put(-5,15){\line(0,1){40}}
\put(145,15){\line(0,1){40}}
\put(10,20){\line(1,0){120}}
\multiput(10,20)(20,0){7}{\line(0,1){20}}
\put(10,41){$\infty$}
\put(30,41){$0$}
\put(50,41){$\xi_1$}
\put(70,41){$\xi_2$}
\put(90,41){$\xi_3$}
\put(110,41){$\xi_4$}
\put(130,41){$\xi_5$}
\thinlines
%$\beta_1$
\put(18,43){$\beta_1$}
\put(5,35){\line(0,1){10}}
\put(5,35){\line(1,0){5}}
\put(25,35){\line(1,0){5}}
\put(5,45){\line(1,0){10}}
\put(15,45){\vector(1,-1){10}}
\put(25,36){1}
\put(35,35){\line(0,1){10}}
\put(10,35){\line(1,0){5}}
\put(30,35){\line(1,0){5}}
\put(25,45){\line(1,0){10}}
\put(25,45){\vector(-1,-1){10}}
\put(13,36){2}

%$\beta_3,\beta_4$
\put(78,43){$\beta_3$}
\put(118,43){$\beta_4$}
\multiput(65,35)(40,0){2}{\line(0,1){10}}
\multiput(65,35)(40,0){2}{\line(1,0){5}}
\multiput(85,35)(40,0){2}{\line(1,0){5}}
\multiput(65,45)(40,0){2}{\line(1,0){10}}
\multiput(85,35)(40,0){2}{\vector(-1,1){10}}
\put(85,36){5}
\put(125,36){5}
\multiput(95,35)(40,0){2}{\line(0,1){10}}
\multiput(70,35)(40,0){2}{\line(1,0){5}}
\multiput(90,35)(40,0){2}{\line(1,0){5}}
\multiput(85,45)(40,0){2}{\line(1,0){10}}
\multiput(75,35)(40,0){2}{\vector(1,1){10}}
\put(73,36){1}
\put(113,36){1}

%\beta_2
\put(38,48){$\beta_2$}
\multiput(0,30)(0,20){2}{\line(1,0){35}}
\multiput(45,30)(0,20){2}{\line(1,0){10}}
\multiput(0,30)(55,0){2}{\line(0,1){20}}
\put(45,30){\vector(-1,2){10}}
\put(35,30){\vector(1,2){10}}
\put(36,30){1}
\put(42,30){5}
\put(20,27){6}
\put(20,51){5}
\put(50,51){1}

%\beta_5
\put(98,48){$\beta_5$}
\multiput(60,30)(0,20){2}{\line(1,0){35}}
\multiput(105,30)(0,20){2}{\line(1,0){35}}
\multiput(60,30)(80,0){2}{\line(0,1){20}}
\put(95,50){\vector(1,-2){10}}
\put(105,50){\vector(-1,-2){10}}
\put(96,30){5}
\put(80,27){3}
\put(80,51){1}
\put(102,30){1}
\put(120,27){3}
\put(120,51){5}
\end{picture}
\caption{$\beta_1,\dots,\beta_5.$}
\label{homb}
\end{figure}

Put 
$$
\begin{cases}
\alpha_1 & = -\sigma_*\rho_* (\beta_1),  \\
\alpha_i & =  \rho_* \beta_i \quad (\text{for}\quad i=2, 3, 4),\\
\alpha_5 & = - \rho_* \beta_5, \\ 
\end{cases}
$$
and $\alpha_i' = \sigma \alpha_i$ and $\beta_i' = \sigma \beta_i$
for $i = 1, \dots , 5$. 
Then we have
$$
\alpha_i \cdot \alpha_j = 0,\quad  \beta_i\cdot \beta_j =0,\quad 
\alpha_i \cdot \beta_j = - \delta_{ij}.
$$
Therefore the 1-cycles $\alpha_1, \dots , \alpha_5, \alpha_1', \dots ,
\alpha_5'$ and $\beta_1, \dots , \beta_5, \beta_1', \dots ,
\beta_5'$ form a symplectic basis for the cup product
on $H_1(C, \bold Z)$.

The inclusion $Prym(C ,\sigma) \to J(C)$ corresponds to the
inclusion $H_1(C, \bold Z)^- \to H_1 (C, \bold Z)$, where
$H_1(C, \bold Z)^-$ is the $(-1)$-part of the action $\sigma_*$.
We have the following proposition.
\begin{proposition}
\label{prop: polarization on prym}
The restriction of
the half of the cup product on 
$H_1(C, \bold Z)$ to $H_1(C, \bold Z)^-$ 
gives a principal polarization on $Prym(C, \sigma)$; that is
$$
A_i \cdot A_j = 0,\quad  B_i\cdot B_j =0,\quad 
A_i \cdot B_j = - 2\delta_{ij},
$$
where
$$
A_i = \alpha_i - \alpha_i',\quad 
B_i = \beta_i - \beta_i'.
$$
\end{proposition}
$H_1(C, \bold Z)^-$ is a free $\bold Z[\rho]$-module of rank 5
generated by $B_1, \dots ,B_5$. 
\begin{remark}
\label{rem:DM list}
By the equation (\ref{eq:degree 6 cov eq}) 
of $C$, $C$ can be regarded as a $\mu_6$
covering of $\bold P^1$ branching at $0, \infty, a_1^2, \dots ,a_5^2$.
The branching index of this covering is 
$\displaystyle (\frac{1}{6},\frac{1}{6},
\frac{1}{3},\frac{1}{3},\frac{1}{3},\frac{1}{3},\frac{1}{3})$
under the notation of \cite{DM},\cite{Mo}.  This index is contained in the
table of \cite{Mo}. This fact gives another proof of the surjectivity
of the period map. (See \S \ref{subsec:three moduli} and 
Theorem \ref{main theorem}.)
Moreover, using the coordinate $a_1^2, \dots ,a_5^2$,
the period map of $J(Y)$ can be expressed by Appell's hypergeometric
functions.
\end{remark}

\subsection{Orthonormal basis and the 27 lines}
\label{subsec:orthonormal bases}

As in \S\ref{subsec:symplectic basis}, we assume that
$0<a_1 < \cdots < a_5$, and regard $C$ as a
6-ple covering of the $\xi$-plane.
Let $\gamma_0$, $\gamma_1$ and $\gamma_i$ $(i=2, \dots 5)$ be the branches of
paths connecting $\infty$ and $0$, $0$ and $\xi_1$, and $\xi_{i-1}$ and
$\xi_i$ $(i=2, \dots ,5)$, respectively,
illustrated in Figure 2.

\setlength{\unitlength}{0.75mm}
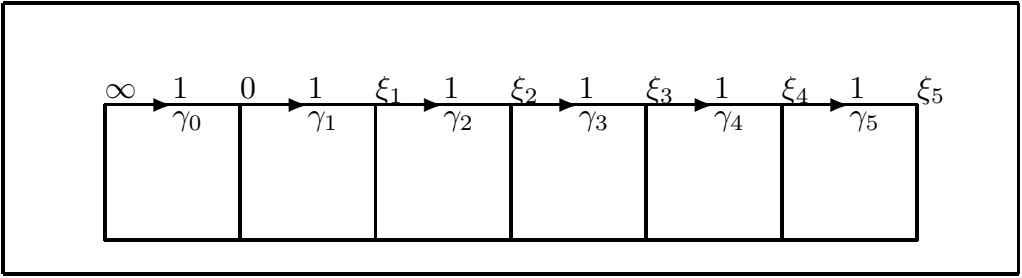
\begin{figure}[hbt]
\begin{picture}(160,50)(-5,10)
%cut lines
\thicklines
\put(-5,15){\line(1,0){150}}
\put(-5,55){\line(1,0){150}}
\put(-5,15){\line(0,1){40}}
\put(145,15){\line(0,1){40}}
\put(10,20){\line(1,0){120}}
\multiput(10,20)(20,0){7}{\line(0,1){20}}
\put(10,41){$\infty$}
\put(30,41){$0$}
\put(50,41){$\xi_1$}
\put(70,41){$\xi_2$}
\put(90,41){$\xi_3$}
\put(110,41){$\xi_4$}
\put(130,41){$\xi_5$}
\multiput(10,40)(20,0){6}{\vector(1,0){10}}
\multiput(20,40)(20,0){6}{\line(1,0){10}}
\put(20,37){$\gamma_0$}
\put(40,37){$\gamma_1$}
\put(60,37){$\gamma_2$}
\put(80,37){$\gamma_3$}
\put(100,37){$\gamma_4$}
\put(120,37){$\gamma_5$}
\put(20,41){$1$}
\put(40,41){$1$}
\put(60,41){$1$}
\put(80,41){$1$}
\put(100,41){$1$}
\put(120,41){$1$}
\end{picture}
\caption{$\gamma_0,\dots,\gamma_5$}
\label{tau}
\end{figure}
An $\bold R$ linear combination of homology classes of a closed paths 
defines an element of the $\bold C$-dual $H^0(C, \Omega^1)^*$
of $H^0(C, \Omega^1)$ by integration and this correspondence 
defines an isomorphism $H_1(C, \bold R) \to H^0(C, \Omega^1)^*$.
Since this isomorphism is equivariant under the action of $\sigma$,
the $(-1)$-eigen space $H_1(C, \bold R)^-$ for $\sigma$ is isomorphic to
the $\bold C$-dual of $H^0(C, \Omega^1)^-$.
A path not necessary closed 
also defines an element of $H^0(C, \Omega^1)^*$. From now on, a path
means the corresponding element in $H_1(C, \bold R)$.
The class of a path from $p$ to $q$ on $C$ in 
$H_1(C,\bold R)/H_1(C,\bold Z)(\simeq J(C))$ corresponds to the divisor class 
$[q]-[p]\in J(C)$ by Abel's theorem. 

The paths $\beta_1, \dots , \beta_5$ can be expressed as
\begin{align*}
&\beta_1 = \gamma_0 - \rho^2\sigma\gamma_0, \\
&\beta_2 = -\rho^2\gamma_0 -\rho^2\gamma_1 +\gamma_1+\rho\sigma\gamma_0, \\
&\beta_3= -\rho^2\gamma_3+\gamma_3, \\
&\beta_4= -\rho^2\gamma_5+\gamma_5, \\
&\beta_5=\gamma_3+\gamma_4+\rho\gamma_5-\rho^2\gamma_5-\rho^2\gamma_4-\rho\gamma_3. 
\end{align*}
By the above equality, we have
$(1-\sigma)\gamma_i \in \frac{1}{1-\rho}H_1(C, \bold Z)^-
=\{ v \in H_1(C, \bold Q)^- \mid (1-\rho )v \in H_1(C, \bold Z)^-\}$.
We put $\bar \gamma_i = (1-\sigma)\gamma_i$ (mod $H_1(C, \bold Z)^-$). 
By the equality
$$
\gamma_0 + \gamma_1 +\gamma_2 +\gamma_3 + \gamma_4 +\gamma_5
-\rho^2\gamma_5 -\rho\gamma_4 -\gamma_3 -\rho^2\gamma_2 -\rho\gamma_1
-\rho^2\sigma\gamma_0 =0,
$$
we have
$$
\bar\gamma_1 = \bar\delta_1+\bar\delta_2, 
\quad \bar\gamma_2 = -2\bar\delta_2+\bar\delta_3-2\bar\delta_5, 
\quad \bar\gamma_3 = \bar\delta_3, 
$$
$$
\bar\gamma_4 = -2\bar\delta_3+\bar\delta_4+\bar\delta_5, 
\quad \bar\gamma_5 = \bar\delta_4, 
$$
where $\bar \delta_i$ is the class of $\frac{1}{1-\rho^2} B_i$
modulo $H_1(C, \bold Z)^-$. 
%%%%%%%%%%%
Since $\sum_{j=1}^i \gamma_j$ is a path from $p_0$ to  $p_i,$ 
$(1-\sigma)\sum_{j=1}^i \gamma_j$ is a path from $\sigma(p_i)$ to $p_i,$ 
which corresponds to the divisor class $[p_i]-[\sigma(p_i)].$
By Proposition {prop:cor vi and ai}, 
we have $v_i=\sum_{j=1}^i\bar\gamma_j$ for
$i=1, \dots ,5$.  More explicitly,
we have
\begin{equation}
\label{eq:base bi and vi}
v_1=\bar\delta_1+\bar\delta_2, \quad
v_2=\bar\delta_1-\bar\delta_2+\bar\delta_3+\bar\delta_5, 
\end{equation}
$$
v_3=\bar\delta_1-\bar\delta_2-\bar\delta_3+\bar\delta_5, \quad
v_4=\bar\delta_1-\bar\delta_2+\bar\delta_4-\bar\delta_5, \quad
v_5=\bar\delta_1-\bar\delta_2-\bar\delta_4-\bar\delta_5.
$$

\begin{definition}
\label{def:hermitian form}
Let $H_3$ be a free $\bold Z[\rho ]$ module equipped with 
a skew symmetric form $\wedge$ 
which satisfies $\rho (v) \wedge \rho (w)= v \wedge w$.
The hermitian metric $h$ on $H_3$ is defined by 
$h(v) =  v \wedge \rho (v)$.
Since the value of 
the associated bilinear form
$q(x,y) = h(x+y)-h(x)-h(y)$ on
$H_3 \times (1-\rho)H_3$ is divisible by 3,
$q$ mod 3 defines a $\bold F_3$ valued bilinear form $q$ on
$H_3/(1-\rho)H_3$.
We denote
$q(\alpha , \alpha )$ by $q(\alpha )$.
Note that $q(\alpha) = 2h(\alpha)$ (mod $3$).
(c.f. \cite{ACT}.)
\end{definition}
We apply this construction to $H_3(Y, \bold Z)$.
Using the basis $A_1, \dots ,B_5$, one can check that
the hermitian form $h$ is isomorphic to 
$$
\left(\begin{matrix}2 & 1 \\ 1 & 2\end{matrix}\right)
^{\oplus 4} \oplus 
\left(\begin{matrix}-2 & -1 \\ -1 & -2\end{matrix}\right)
$$ 
as a quadratic form.
Moreover, $H_3(Y, \bold Z)/(1-\rho)H_3(Y, \bold Z)\simeq \bold F_3 ^5$
is equipped with a $\bold F_3$ valued quadratic form $q$.
The multiplication by $1-\rho^2$ induces an isomorphism from
$\frac{1}{1-\rho}H_3(Y, \bold Z)/H_3(Y, \bold Z)$ to 
$H_3(Y, \bold Z)\otimes \bold Z[\rho]/(1-\rho)$
and via this isomorphism $q$ is considered as a quadratic form on 
$$
\frac{1}{1-\rho}H_3(Y, \bold Z)/H_3(Y, \bold Z)
\simeq J(Y)_{1-\rho},
$$
where $J(Y)_{1-\rho}$ is
the $(1-\rho)$-torsion subgroup of $J(Y)$.
Let $\bar \delta_i$ be the image of $\frac{1}{1-\rho^2} B_i$
in $\frac{1}{1-\rho}H_3(Y, \bold Z)/H_3(Y, \bold Z)$ under the map 
$\bar c$ in Theorem \ref{theorem:prym}.
By the compatibility of the symplectic forms on 
$H_1(C, \bold Z)^-$ and $H_3(Y, \bold Z)$
stated in Remark \ref{rem:pol J and C},
the quadratic form $q$ is
$$
q(\sum_{i=1}^5 \bar t_i \bar \delta_i) = 
-\sum_{i=1}^4 \bar t_i^2 + \bar t_5^2.
$$
By using this explicit formula and (\ref{eq:base bi and vi})
we can check that $\{ v_1, \dots , v_5 \}$ is an orthonormal basis.

We interpret the geometry of the 27 lines in a cubic surface into
the finite geometry over $\bold F_3^5$ with the quadratic form $q$.
Let $\{ L_1, L_2,L_3\}$ be a tritangent of the 27 lines in $X$, i.e.
$L_1$, $L_2$ and $L_3$ are contained in a hyperplane in $\bold P^3$.
By choosing a marking of $X$, we assume that $L_1 = E_1$, $L_2 = L_{16}$ 
and $L_3= C_6$.  Then we have $v_1=\Lambda ([L_1]-[L_2])$. 
Since $[C_6]$ is fixed
under the Clemens-Griffiths involution and it acts as $(-1)$-multiplication
on the image of $\Lambda$, we have 
\begin{equation}
\label{eq:tritang and level map}
\Lambda ([E_1]-[C_6])=-\Lambda ([L_{16}]-[C_6]).
\end{equation}
Thus $\Lambda (2[C_6]-[L_{16}]-[E_1])=0$.  It is easy to see that
$[E_i]-[L_{i6}]$ $(i=1, \dots , 5)$ and $2[C_6]-[L_{16}]-[E_1]$ generates 
$H_2^{prim}(X, \bold Z)\otimes (\bold Z/3\bold Z)$ 
by direct calculation.  Using this basis
of $H_2^{prim}(X, \bold Z)\otimes (\bold Z/3\bold Z)$,
we can show that $q(\Lambda (v))= (v,v)$ (mod $3$),
where $(,)$ denotes the intersection form on $H_2^{prim}(X, \bold Z)$.
Therefore by
the map $\Lambda$, $J(Y)_{1-\rho}$ is identified with the quotient of
$H_2^{prim}(X, \bold Z)\otimes (\bold Z/3\bold Z)$ by the radical 
of the intersection form mod 3 on $H_2^{prim}(X, \bold Z)$.
Since the action of $W(E_6)$ on $H_2^{prim}(X, \bold Z)$ preserves
the intersection form, this action induces a linear map 
on $J(Y)_{1-\rho}$ preserving
the quadratic form $q$. Thus we have a map $W(E_6) \to O(5,\bold F_3)$,
where $O(5, \bold F_3)=\{ g \in GL (5, \bold F_3) \mid g ^tg = I\}$.
The composite
$Aut(\Gamma_{std}) \to O(5, \bold F_3) \to 
PO(5, \bold F_3)=O(5, \bold F_3)/\{\pm I\}$ is an isomorphism.

The inverse of the above map is explicitly given as follows.
By the equality (\ref{eq:tritang and level map}), we have
$$
\pm \Lambda ([E_1]-[L_{16}])=\pm \Lambda ([E_1]-[C_6])
= \pm \Lambda ([L_{16}]-[C_6]).
$$
This equality implies that the map $t$ from the set of tritangents to
the subset
$$
T=\{ v = ^t(u_1, \dots ,u_5)\in \bold F_3^5\mid 
q(v) = 1 \}/\{\pm 1\}
$$
of $\bold F_3^5/\{\pm 1\}$ is well defined.  
Two elements $v_1$ and $v_2$ in $T$ are said to be vertical if
and only if $q(v_1, v_2)=0$.
\begin{proposition}
\label{prop:bijection}
The map 
$$
t : \{ \text{ tritangents of $X$}\} \to T
$$
is bijective. Moreover by this bijection, 
two tritangents are collinear, i.e. there exists a common line
in them, if and only if the corresponding elements in
$T$ are vertical. 
\end{proposition}
\begin{proof}
Since $\#\{ \text{ tritangents of $X$}\}=\# T =45$,
it is enough to
prove the injectivity of the map $t$.  If $T_1=\{ N_1, N_2,N_3\}$
and $T_2=\{M_1, M_2,M_3\}$ are collinear, then we showed that 
$t(T_1)$ and $t(T_2)$ are vertical to each other, therefore the images
are different.  If $T_1$ and $T_2$ are not collinear, then
we may assume that $N_i \cap M_i \neq \emptyset$ for $i=1,2,3$ and
$N_i \cap M_j = \emptyset$ if $i \neq j$.  Since two lines
$N_2$ and $M_1$ are disjoint, there are exactly 5 lines which
intersect both $N_2$ and $M_1$.  Two of them are $N_1$ and $M_2$.
The rest of them are written as $F_3$, $E_4$, $E_5$.
Then we can find $E_3$ such that $E_3, F_3, N_2$ is a 
tritangent. We denote $N_1= E_1$, $M_2= E_2$.
Then $E_1, \cdots, E_5$
are mutually disjoint lines.  Moreover, there exists a unique line 
$E_6$ such that
$E_1, \dots ,E_6$ are mutually disjoint lines.  By blowing them
down, we obtain a marking of $X$.
In particular, $N_2=C_6$, $M_1=C_3$,
$T_1 =\{ E_1, C_6,L_{16}\}$ and $T_2=\{ E_2,C_3, L_{23}\}$.
Note that the action of $\frak S_5$ on the set of tritangents and
$T$ are equivariant. The action of the transposition
$(14)$ is trivial on $t(T_2)$ and non-trivial on $t(T_1)$, therefore
$t(T_1)$ and $t(T_2)$ are different.
\end{proof}

The explicit correspondence
between the set of tritangents and the set $T$ is given as
\begin{align*}
& t(C_6E_iL_{i6})=v_i,\quad
t(E_6C_iL_{i6})=\sum_{j\neq i}v_j, \\
& t(E_iC_jL_{ij})=v_j-\sum_{k\neq i,j}v_k, \quad (1 \leq i,j \leq 5), \\
& t(L_{ij}L_{kl}L_{m6})=v_i+v_j-v_k-v_l. 
\end{align*}
By this correspondence,
the group $PO(5, \bold F_3)$ acts on this set $T$.

Since the set of the 27 lines in a cubic surface can be identified with 
$$
\Cal L_{cs} = \{ \mu =\{ T_1, \dots , T_5\} \mid 
\begin{minipage}{5.2cm}
$T_i$ are distinct tritangents\\
\noindent
which are colinear to each other.
\end{minipage}
\}, 
$$
it is also identified with the set $\Cal L$
$$
\Cal L = \{ \nu =\{ T_1, \dots , T_5\} \subset T \mid \text{  
$T_i$ is vertical to $T_j$ for $1 \leq i < j \leq 5$} \}.
$$
The two lines $\mu$ and $\mu'$ corresponding to two elements 
$\nu$ and $\nu'$ in $\Cal L$
intersect if and only if $\nu$ and
$\nu'$ have a common element of $T$. 
Via this identification, $PO(5, \bold F_3)$ acts on the graph
$\Gamma_{std}$ and we have the map $PO(5, \bold F_3) \to Aut(\Gamma_{std})$.
This map is the inverse of the map $Aut(\Gamma_{std}) \to PO(5, \bold F_3)$.

For $1\leq i<j \leq 6$, the transposition of the index $i,j$
for $e_i, c_i, l_{ij}$ in $\Gamma_{std}$
induces an involution of $\Gamma_{std}$,
which is denoted by $r_{ij}$. 
The involutions $r_{ij}$
generate the symmetric group $\frak S_6$.  The element $r_{123}$
of $Aut(\Gamma_{std})$ determined by
\begin{align*}
& r_{123}(E_i)=L_{jk}, r_{123}(C_i)=C_{i} 
\text{ for }\{ i,j,k\}=\{ 1,2,3\}, \\
& r_{123}(C_i)=L_{jk}, r_{123}(E_i)=E_{i} 
\text{ for }\{ i,j,k\}=\{ 4,5,6\}, \\
& r_{123}(L_{ij})=L_{ij} \text{ for }i \in \{ 1,2,3\}
\text{ and } j\in\{ 4,5,6\} 
\end{align*}
is an involution.
Involutions $r_{ij}$ and $r_{123}$ generate $Aut(\Gamma_{std})$.

On the other hand, for $1\leq i < j \leq 5$,
the linear map $R_{ij}\in PO(5, \bold F_3)$
given by 
$$
R_{ij}(v_k) = v_k \ (i,j \neq k), \quad
R_{ij}(v_i) = v_j, \quad
R_{ij}(v_j) = v_i
$$
is an element of $O(5, \bold F_3)$. In general, for a vector $v$
such that $q(v)=2$, we define a reflection $R_v \in O(5, \bold F_3)$ 
with respect to $v$
by $R_v\mid_{v\bold F_3}=-1$ and $R_v\mid_{v^{\perp}}=1$.
The reflections $R_v$ for $v=(1,1,1,1,-1)$ and $(1,1,1,-1,-1)$
are denoted by $R_{56}$ and $R_{123}$, respectively.

\begin{lemma}
\label{compatibility}
Under the isomorphism $\iota: Aut(\Gamma_{std}) \to PO(5, \bold F_3)$, 
the involution $r_{ij}$ 
in $Aut(\Gamma_{std} )$ corresponds to the 
involution $R_{ij}$ in $PO(\bold F_3, 5)$
for $1 \leq i < j \leq 5$. Moreover the involutions $r_{56}$ 
and $r_{123}$ correspond
to $R_{56}$ and $R_{123}$ respectively. 
\end{lemma}
\begin{proof}
The automorphism of $H_2^{prim}(X, \bold Z)$ induced by $R_{56}$ (resp. $R_{123}$)
is the reflection with respect to the element $[E_5]-[E_6]$ (resp. 
$[H]-[E_1]-[E_2]-[E_3]$).
By the definition of the isomorphism between $Aut(\Gamma_{std})$ and
$PO(5, \bold F_3)$, $r_{56}$ (resp. $r_{123}$) is the reflection with respect to
$\Lambda ([E_5]-[E_6])$ (resp. $\Lambda ([H]-[E_1]-[E_2]-[E_3])$). 
Using the compatibility of $q$ and the intersection
form on $H_2^{prim}(X, \bold Z)$, we get the reflection vector for $r_{56}$
(resp. $r_{123}$).
\end{proof}

\section{Zero of Theta functions restricted to curves}
\label{sec:zero of theta}

\subsection{Moduli space of abelian varieties with $\mu_3$ actions}
\label{subsec:moduli of abel var with rho action}
We consider a moduli space of a certain kind of abelian
varieties and its analytic expression.
Let $J$ be an abelian variety. An element of $H^2(J, \bold Z)$ is called
a polarization if it is the Chern class of an ample divisor. By 
the isomorphism $H^2(J, \bold Z) \simeq 
Hom(\wedge^2 H_1(J, \bold Z) \to \bold Z)$, it corresponds to a symplectic
form on $H_1(J, \bold Z)$.  For a principally polarized abelian variety $J$
with an action of $\mu_3$,
the representation of $\mu_3$ on $H^0(\Omega^1)$ is 
called the type of the action of $\mu_3$.
Let us denoted by $\chi$ the natural representation of $\mu_3$ 
on $\bold C$.  
For a principally polarized abelian variety with an action of $\mu_3$
of type $4\chi \oplus \bar\chi$,
a natural $\bold F_3$-valued
quadratic form $q$ is defined
on the $(1-\rho)$-torsion part $J_{1-\rho}$ of 
$J$ as in Definition \ref{def:hermitian form}.
A level $(1-\rho )$ structure is defined as an isomorphism
$\Psi_{ab} :\oplus_{i=1}^5 \bold F_3v_i \to J_{1-\rho}$ 
such that $\{ \Psi_{ab} (v_1), \dots ,\Psi_{ab} (v_5)\} $ is
an orthonormal basis.
\begin{definition}
\label{def:moduli ab}
The set of isomorphism classes of triples 
$(J, \iota , \Psi_{ab})$ is denoted by $\Cal M_{ab}$,
where $J$ is a
principally polarized abelian variety, $\iota$ is a $\mu_3$ action of type
$4\chi \oplus \bar\chi$ and $\Psi_{ab}$ is a level $(1-\rho)$ structure.
\end{definition}
\begin{remark}
\label{def:o53 action}
The orthogonal group $O(5, \bold F_3)$ acts as the right multiplication
of $\oplus_{i=1}^5\bold F_3 v_i$. 
Through the equality $g\Psi_{ab} (v)= \Psi_{ab} (vg)$,
we define the left action of $O(5,\bold F_3)$ on $\Cal M_{ab}$.
Since 
$(J, \iota , \Psi_{ab})$ is isomorphic to 
$(J, \iota , -\Psi_{ab})$ as triples, this action reduces
to the action of $PO(5, \bold F_3)$ which is isomorphic to $W(E_6)$. 
\end{remark}
To obtain an analytic expression of $\Cal M_{ab}$, we introduce the notion of
homology marking.
We define the standard module $Mod_{std}$ as
the free $\bold Z$ module generated by $a_1, \dots , a_5,b_1, \dots ,b_5$.
On the module $Mod_{std}$, we define a symplectic form $\phi$ as
$\phi (a_i, a_j ) =\phi (b_i, b_j ) = 0$ and $\phi (a_i, b_j )=-\delta_{ij}$
and the action of $\rho$ as
$$
\left(
 \alpha_1 , \dots  , \alpha_5 , \beta_1 , \dots , \beta_5 
\right)
\left(
\begin{matrix} a_1  \\  \vdots   \\  a_5 
\\  b_1 \\  \vdots \\ b_5 \\
\end{matrix}
\right)
\mapsto
\left(
 \alpha_1 , \dots  , \alpha_5 , \beta_1 , \dots , \beta_5 
\right)
 W
\left(
\begin{matrix} a_1  \\  \vdots   \\  a_5 
\\  b_1 \\  \vdots \\ b_5 \\
\end{matrix}
\right),
$$
where 
$$
W=\left(
\begin{matrix}
-I &- H \\
H & 0 \\
\end{matrix}
\right), \quad
H = 
diag (1,\dots ,1,-1).
$$

For a principally polarized abelian variety $J$ with an action of $\mu_3$ of
type $4\chi\oplus \bar\chi$, a homology marking is defined
as an isomorphism $\Psi_{hom} : Mod_{std} \to H_1(J, \bold Z)$
compatible with the symplectic forms and the action of $\rho$.
The triple $(J, \iota , \Psi_{hom} )$ is called the homology marked abelian
variety with an action of $\mu_3$ of type $4\chi \oplus \bar\chi$.
\begin{definition}
\label{def:moduli ab hom}
The set of isomorphism classes of
homology marked abelian varieties
$(J, \iota , \Psi_{hom})$ is called the moduli
spaces of homology marked abelian variety and
denoted by $\Cal M_{hom}$.
\end{definition}

The opposite automorphism group
$Aut^0(Mod_{std})$ is defined as the copy of  $Aut(Mod_{std})$
whose product structure is defined by the reverse of that of $Aut(Mod_{std})$.
The group $Aut^0(Mod_{std})$
acts on $Mod_{std}$ from the right.
The group $Aut^0(Mod_{std})$ acts on $\Psi_{hom}$ from the left
by the rule $(g\Psi_{hom})(v) = \Psi_{hom} (v g)$.  
If we use the basis $a_1, \dots , b_5$,
$Mod_{std}$ can be identified with the 10 dimensional integer valued row vector
and $Aut^0(Mod_{std})$ is identified with the centralizer of $W$ 
in $Sp(10, \bold Z)$. 
An element  $g$ in $Aut^0(Mod_{std})$ acts on $\Cal M_{hom}$ by  
$g(J,\iota,\Psi_{hom})=(J,\iota,g\Psi_{hom})$.

For a homology marking $\Psi_{hom} : Mod_{std} \to H_1(J, \bold Z)$,
$\{ A_1=\Psi_{hom} (a_1), \dots , B_5= \Psi_{hom} (b_5)\}$ becomes a 
symplectic basis.
By using this basis, the module $H_1(J, \bold Z)$ 
is identified with the integer valued row vector space $\bold Z^{10}$.
Using these data, 
we obtain a period matrix $\tau$ in 
$\frak H_5=\{ \tau \in M(5, \bold C) \mid \tau =\;^t\tau, \Im (\tau )
\text{ is positive definite }\}$ such that 
$W \cdot \tau =\tau$ as follows.
Let $\phi_i$ be a normalized 1-form on $J$. i.e. 
$\displaystyle \int_{B_i}\phi_j = \delta_{ij}$ and set
$\displaystyle \tau_{ij}= \int_{A_i} \phi_j$.  In other words,
$$
\left(
\begin{matrix}
A_1 \\
\vdots \\
A_5 \\
B_1 \\
\vdots \\
B_5 \\
\end{matrix}
\right)
\left(
\begin{matrix}
\phi_1 &,\dots, & \phi_5 \\
\end{matrix}
\right)
=
\left(
\begin{matrix}
\tau \\
I \\
\end{matrix}
\right).
$$
Here we write 
$A_i\cdot \phi_j = \displaystyle \int_{A_i} \phi_j$
for simplicity. 
The matrix $(\tau_{ij})_{ij}$ is called the normalized period matrix.
Then the normalized period matrix
$\tau$ is symmetric and $\Im (\tau )$ is 
positive definite.
Let $S$ be the matrix for the $\rho^*$ action on $(\phi_i)_i$, i.e.
$$
\rho^*(\phi_1, \dots ,\phi_5) =(\phi_1, \dots ,\phi_5)S.
$$
Then we have the following relation:
\begin{align*}
\left(
\begin{matrix}
\tau \\
I \\
\end{matrix}
\right)
= &
\rho_*
\left(
\begin{matrix}
A_1 \\
\vdots \\
A_5 \\
B_1 \\
\vdots \\
B_5 \\
\end{matrix}
\right)
(\rho^{-1})^*
\left(
\begin{matrix}
\phi_1 &, \dots ,  & \phi_5 \\
\end{matrix}
\right)  
= 
W
\left(
\begin{matrix}
A_1 \\
\vdots \\
A_5 \\
B_1 \\
\vdots \\
B_5 \\
\end{matrix}
\right)
\left(
\begin{matrix}
\phi_1 & ,\dots,  & \phi_5 \\
\end{matrix}
\right)   S^{-1} \\
= &
W \left(
\begin{matrix}
\tau \\
I \\
\end{matrix}
\right) S^{-1}.
\end{align*}
Therefore we have $S= H\tau$ and $W \cdot \tau = \tau$, i.e.
$(H\tau )^2 + H\tau +I = 0$.
Here an element 
$$
g=\left(\begin{matrix}
A & B \\ C & D \\
\end{matrix}\right)
$$
in $Sp(10, \bold Z)$ acts on $\frak H_5$ by
$g\cdot \tau = (A \tau +B)(C \tau + D)^{-1}$ for $\tau \in \frak H_5$.
Thus on the period lattice 
$L  = \bold Z^5 \tau \oplus \bold Z^5 \subset \bold C^5$, the action
of $\rho$ is given by the right multiplication of $S$.
Under these notation, the action of $\rho$ on $J=
\bold C^5/ (\bold Z^5\tau \oplus \bold Z^5)$ is also given by
the right multiplication of $S$.

\subsection{Moduli space of abelian varieties as a ball quotient}
\label{subsec:ball quotient}
We give an isomorphism between 
the moduli space $\Cal M_{hom}$ and 4-dimensional complex ball 
$\bold B_4= \{^t(x_1, \dots ,x_4) \mid 
\mid x_1 \mid^2 +\cdots +\mid x_4 \mid^2
<1 \}$.

Since the matrix $H\tau$ represents $\rho^*$
on $H^0(J, \Omega^1)$ with respect to the normalized 1-forms,
there exists a unique eigen vector 
$\eta =^t(\eta_1,\cdots ,\eta_5)  \in \bold C^5$ 
up to constant corresponding to the eigen value
$\bar\omega$ of the 
left multiplication of $H\tau$.
In other words, we have $ H\tau\eta =\bar\omega \eta$.
Since $\Im (\tau ) > 0$, we have
\begin{align*}
0 & < ^t\bar\eta  \Im (\tau ){\eta} =
{^t\bar\eta} \frac{\tau-\bar\tau}{2i} {\eta} \\
& =\frac{1}{2i}(-\omega{^t\bar\eta H} \eta +
\omega^2{^t\bar\eta H} \eta) \\
& = -\frac{\sqrt 3}{2}({^t\bar\eta} H\eta). 
\end{align*}
Thus we have 
$\mid \eta_1 \mid^2 +\cdots + \mid \eta_4 \mid^2 - \mid \eta_5 \mid^2 <0$,
so $x =\; ^t(x_1, \dots, x_4)$ is 
an element of $\bold B_4$, where  
$\displaystyle x_i = \frac{\eta_i}{\eta_5}$
($i=1, \dots , 4$).
We can regard $\bold B_4$ as an open set of 
$\bold P^4=\{\;^t(\eta_1:\cdots :\eta_5)\}$.
Conversely we obtain the normalized period matrix
$\tau$ from the vector $x$ in $\bold B_4$. For any vector
$\delta \in \bold C^5$ in the $\omega$-eigen space for the 
left multiplication of $H\tau$, we have
$^t\delta H {\eta}=0$ because of the equalities
$$
\omega \;^t\delta H \eta = \;
^t\delta \tau \eta = \bar\omega \;^t\delta H \eta.
$$
Therefore the matrix $H\tau$ can be
characterized as 
$$
H\tau \mid_{\eta\bold C}= \omega\cdot I,\quad 
H\tau \mid_{\eta^{\perp_H}} = \bar\omega\cdot I,
$$
where $\eta^{\perp_H}$ is the vector space
$\{ \delta \mid \; ^t\delta H {\eta}=0\}$.
As a consequence, we have the isomorphism
$$
\bold B_4 \simeq \{ \tau \in \frak H_5 \mid W\cdot \tau = \tau \}
(\simeq \Cal M_{hom}).
$$

The module $Mod_{std}$ is freely generated 
by $B_1, \dots ,B_5$ over $\bold Z[\rho]$. Since 
the action of $Aut^0(Mod_{std})$ is 
$\bold Z[\rho]$-linear, $Aut^0(Mod_{std})$ can be identified with
a subgroup $U(4,1)$ of $GL(5, \bold Z[\rho])$ using the $\bold Z[\rho]$
basis $B_1, \dots , B_5$. 
Here $U(4,1)=\{ g \in GL(5,\bold Z[\rho ]) \mid\; ^t\bar g H g = H \}$.
For an element
$$
g= A +B \rho \quad (A, B \in M(5, \bold Z))
$$
of $U(4,1)$, the corresponding element $\iota (g)$ 
in $Sp(10,\bold Z)$ is given by 
\begin{equation}
\label{eq:group embedding}
\iota (g) =
\left(
\begin{matrix}
H(A-B)H & -HB \\
BH & A \\
\end{matrix}
\right).
\end{equation}
We define the orthogonal group $O(4,1)$ by $U(4,1) \cap M(5, \bold Z)$.
Note that as in $\Cal M_{ab}$,
the action of $-\rho I \in U(4,1)$ on $\Cal M_{hom}$ is trivial and
the action of $U(4,1)$ on $\bold B_4$ factors through the action of 
$PU(4,1)=U(4,1)/ < -\rho I>$.

A homology marking $\Psi_{hom}:Mod_{std} \to H_1(J, \bold Z)$, 
induces an isomorphism
$$
\frac{1}{1-\rho}Mod_{std}/Mod_{std} \simeq J_{1-\rho}.
$$
Therefore the elements $v_1, \dots ,v_5$ in $\frac{1}{1-\rho}Mod_{std}/Mod_{std}$
given by the equality 
(\ref{eq:base bi and vi})
form an orthonormal basis in $\frac{1}{1-\rho}Mod_{std}$
with respect to the $\bold F_3$-valued quadratic form $q$. 
Using this basis, we get an isomorphism
$\Psi_{ab} : \oplus \bold F_3 v_i \to J_{1-\rho}$
preserving the $\bold F_3$-valued quadratic form.  This correspondence
gives a morphism of moduli spaces $\Cal M_{hom} \to \Cal M_{ab}$.
We define the congruence subgroup $\Gamma (1-\rho )$ by
$$
\Gamma (1-\rho) = \{ g \in U(4,1)\subset GL(5, \bold Z[\rho])
\mid g \equiv I \text{ (mod $(1-\rho)$)}\}.
$$
By the identification $\Cal M_{hom} \simeq \bold B_4$,
$\Cal M_{ab}$ is identified with the quotient of $\bold B_4$ by the congruence
subgroup $\Gamma (1-\rho)$ of $U(4, 1)$.
The action of $A + B \rho \in U(4,1)$ on $\bold B_4 \subset \bold P^4$
is given by $^t(\eta_1, \dots ,\eta_5)\mapsto 
(A+B\omega^2 )\;^t(\eta_1, \dots ,\eta_5)$.

\subsection{Theta functions and their transformation formulae}
\label{subsec:transformation formula}

We recall the transformation formula for the acton of $Sp(10,\bold Z)$ 
on theta constants (see \cite{I}).
\begin{definition}
\label{def:theta const}
For $(m',m'') \in \bold R^{10}$ and $\tau \in \frak H_{5}$,
we define theta functions $\Theta_{m}(\tau, z)$, $\Theta (\tau, z)$
and theta constants $\Theta_m(\tau )$ as
\begin{align*}
& \Theta_{(m', m'')}(\tau, z) = \sum_{p \in \bold Z^5}
\bold e (\frac{1}{2}(p+m')\tau ^t(p+m') + (p+m')^t(z+m'')), \\
& \quad \Theta (\tau, z ) = \Theta_{(0,0)}(\tau,z), \\
& \Theta_{(m',m'')}(\tau) = \Theta_{(m',m'')}(\tau,0). 
\end{align*}
\end{definition}
\begin{proposition}
\label{prop:theta functional equation}
For any 
$$
g = \left(
\begin{matrix}
A & B \\
C & D \\
\end{matrix}
\right) \in Sp(10, \bold Z)
$$
and $m = (m' , m'') \in\bold R^{10}$, put
\begin{align*}
(\tau^\#,z^\#) = & ((A\tau + B)(C\tau +D)^{-1},z(C\tau +D)^{-1}), \\
m^\# = & m g^{-1} + \frac{1}{2}((C{^tD})_0, (A{^tB})_0). 
\end{align*}
Then we have
\begin{equation}
\label{eq:theta trans}
\Theta_{m^\#}(\tau^\#, z^\#) =
\bold e (\frac{1}{2}z(C\tau + D)^{-1}C{^tz})
\det (C \tau + D)^{1/2}
u(g) \cdot \Theta_{m}(\tau,z )
\end{equation}
where $u(g) \in \bold C^{\times}$ depends only on $g$.
\end{proposition}
\begin{proposition}
\label{prop:rho transform}
Let $\tau \in \frak H_5$ be the normalized period matrix of 
an abelian variety with an action of $\mu_3$ of type 
$4\chi \oplus \bar\chi$.
Then we have
$$
\Theta_{\frac{1}{2}\bold 1,\frac{1}{2}\bold 1}(\tau , z(H\tau)^{-1})
=\bold e (\frac{1}{2}z \tau^{-1}z)
\Theta_{\frac{1}{2}\bold 1,\frac{1}{2}\bold 1}(\tau , z).
$$
\end{proposition}
\begin{proof}
Since $H\tau$ is equal to the representation matrix of
$\rho^*$ for the normalized 1 forms
$\phi_1, \dots ,\phi_5 \in H^0(J, \Omega^1)$,
we have $\det (H\tau)=1$. 
By applying the transformation formula in
Proposition \ref{prop:theta functional equation}
for $g = W$ in \S \ref{subsec:moduli of abel var with rho action},
we have
\begin{align*}
\Theta_{\frac{1}{2}\bold 1,\frac{1}{2}\bold 1}(\tau , z(H\tau)^{-1})
= & 
\bold e (\frac{1}{2}z\tau^{-1}{^tz})u(W)
\Theta_{-\frac{1}{2}d(H), d(H)-\frac{1}{2}\bold 1}(\tau , z) \\
= & 
\bold e (\frac{1}{2}z\tau^{-1}{^tz})u(W)
\Theta_{\frac{1}{2}\bold 1, \frac{1}{2}\bold 1}(\tau , z). 
\end{align*}
By evaluating this equality at
$\tau = diag (\omega ,\cdots ,\omega , -\omega^2)$,
and $z = \bold 1 ((H\tau)^{-1}-I)^{-1}$, we get
$u(W)=1$. Thus we get the proposition.
\end{proof}

\begin{definition}
\label{def:theta div}
Let us fix a period matrix $\tau \in \frak H_5$.
The zero $\{ z \in \bold C^5 \mid \Theta (z, \tau)=0\}$ of $\Theta$ is
invariant under the translation of the period lattice 
$\bold Z^5 \tau \oplus \bold Z^5$.  Therefore zero locus of $\Theta$
determines a divisor $D_{Prym, \Theta}$ 
of $\bold C^5/(\bold Z^5 \tau \oplus \bold Z^5)$.  
It is called the theta divisor for this polarization.  
\end{definition}

As a corollary of Proposition \ref{prop:rho transform},
we have the following proposition.
The vector $\bold 1\tau + \bold 1$ is denoted by $\newfont I$.
\begin{proposition}
\label{prop:invar theta div for rho}
Let $\tau$ be an element in $\frak H_5$ satisfying the condition of
Proposition \ref{prop:rho transform}.
Then the divisor $D_{Prym, \Theta}+\frac{1}{2}\newfont I$
is stable under the action of $\rho$.
\end{proposition}

\subsection{Zero of the restricted theta function}
\label{subsec:restricted theta}

In this subsection, we assume that X has no Eckardt points
and apply the results in \S \ref{subsec:moduli of abel var with rho action}
\ref{subsec:ball quotient}, \ref{subsec:transformation formula}
to the Prym variety $Prym(C, \sigma)$ for the double coverings of \S 
\ref{subsec:cylinder map}.
Let $C$ be a curve defined by the equation \ref{eq:degree 6 cov eq}
with $a_i \in \bold R$ for $i=1, \dots , 5$ and $0< a_1 < \cdots < a_5$.
We consider the symplectic basis 
$A_1, \dots ,A_5,B_1, \dots , B_5$ 
for the half of the cup product of $H_1(C, \bold Z)^-$
defined in 
\S \ref{subsec:symplectic basis}.
Let $\phi_1, \dots ,\phi_5$ be a basis of the $(-1)$-eigen space of 
holomorphic 1-form for the action of $\sigma$.
This is said to be normalized if 
$$
\int_{B_i} \phi_j = \delta_{ij}.
$$
Then $A_1, \dots , B_5$ and $\phi_1, \dots ,\phi_5$ are regarded as
a symplectic basis for $H_1(Prym(C, \sigma), \bold Z)$ and
a normalized 1-form on $Prym(C, \sigma)$ 
for the symplectic base $A_1, \dots ,B_5$.
We define the normalized period matrix $\tau$ as in the last subsection.
Then the Prym variety $Prym(C, \sigma)$
is isomorphic to the complex torus 
$\bold C^5/(\bold Z^5 \tau \oplus \bold Z^5)$.

Let $\gamma$ be a path in $C$ connecting $Q_0$ and $p$.
The integral
of $\phi=(\phi_1, \dots , \phi_5)$ along this path is 
written by $\int_{Q_0}^x\phi$ for short.
Since $Q_0$ is fixed under the action of $\sigma$, the path 
$\gamma -\sigma (\gamma )$ connects $\sigma (p)$ and $p$.
The integral $\int_{\sigma (p)}^p \phi$ 
of $\phi$ along $\gamma -\sigma (\gamma )$ is a $\bold C^5$ valued
holomorphic function on the universal covering of $C$.
The theta function $\Theta_m(\tau , z)$ is denoted
by $\Theta_m(z)$ for simplicity.
For $v \in \bold C^5$, we consider a holomorphic function
$\Theta_{\frac{1}{2}\bold 1, \frac{1}{2}\bold 1} (v + \int_{\sigma (p)}^p \phi)$ 
on the universal covering of
$C$.  By the quasi periodicity of theta functions, the order
of zero of the function 
$\Theta_{\frac{1}{2}\bold 1,\frac{1}{2}\bold 1} (v + \int_{\sigma (p)}^p \phi)$
depends only on $p$ in $C$.

First we investigate
the order of zero of the function 
$\Theta_{\frac{1}{2}\bold 1,\frac{1}{2}\bold 1} (v + \int_{\sigma (p)}^p \phi)$
at $p\in \Sigma = \{ p_1, \dots ,p_5,\sigma (p_1), \dots ,\sigma (p_5),
p_0, p_{\infty} \}$,
where 
$v \in \frac{1}{1-\rho}L=
\frac{1}{1-\rho}(\bold Z^5 \tau \oplus \bold Z^5)$. 
The class of $\int_{\sigma (p)}^p\phi$ for $p \in \Sigma$
in $\frac{1}{1-\rho}L/L$ is denoted by $\jmath (p)$.
Then the map $\jmath: C \to \bold C^5/L \simeq Prym(C, \sigma)\simeq J(Y)$
is equal to the map $\jmath$ defined in 
Proposition \ref{prop:cor vi and ai}. Therefore
$\jmath (p) \in \frac{1}{1-\rho}L/L$ for $p \in \Sigma$ and
$\jmath (p_i)=v_i$ and $\jmath (\sigma (p_i))=-v_i$.
Using the $\bold F_3$ valued quadratic form $q$
introduced in Definition \ref{def:hermitian form},
the order of zero of theta function modulo 3 is
give by the following proposition.
\begin{proposition}
\label{prop:zero order mod 3 of theta}
The order of zero of $\Theta(\frac{1}{2}\newfont I +
v+\int_{\sigma (p)}^p\phi)$ at $p \in \Sigma$
is equal to $q(v+\jmath (p))$ modulo 3.
\end{proposition}

Next, we study the zeros of the pull back of the theta functions
by the $\mu_3$ equivariant maps $\jmath +v$ for 
$v \in Prym(C, \sigma)$. 
Even though the theta divisor depends on the choice of
symplectic basis of $H_1(Prym(C, \sigma), \bold Z)$, 
the cohomology class
$c_1(D_{Prym, \Theta}) \in H^2(Prym(C, \sigma), \bold Z(1))$
of the theta divisor is independent of the choice of symplectic basis
and it is equal to 
$\sum_{i=1}^5A_i^* \wedge B_i^*$, where $\{ A_1^*, \dots , A_5^*, 
B_1^*, \dots ,B_5^* \}$ is the dual basis of the symplectic basis.
The degree of the inverse image of $D_{Prym ,\Theta}$
by $\jmath$ is equal to the image of $c_1(D_{Prym, \Theta})$ under the morphism
$$
H^2(Prym(C, \sigma), \bold Z(1)) \overset{\jmath^*}\longrightarrow
H^2(C, \bold Z(1)) \overset{\deg}\longrightarrow \bold Z.
$$
Using the relation
\begin{align*}
\jmath^*(A_i^*)=& \alpha_i^*-{\alpha_i'}^*, \\
\jmath^*(B_i^*)=& \beta_i^*-{\beta_i'}^*, 
\end{align*}
and
$\deg (jac^*(\alpha_i^*\wedge \beta_i^*))=1 $ for $jac:C\to J(C)$ defined in 
\S\ref{subsec:abel-Jacobi and orthogonal} and  $i=1, \dots ,5$,
we have  $\deg (\jmath^*(\sum_{i=1}^5A_i^* \wedge B_i^*))= 10$.
As a consequence, we have the following proposition.
\begin{proposition}
\label{prop:degree is 10}
Let $v \in Prym(C, \sigma)$. Suppose that 
the image of $\jmath + v: C \to Prym(C, \sigma)$
is not contained in the theta divisor $D_{Prym, \Theta}$.
Then the degree of the inverse image of $D_{Prym, \Theta}$
by the map $\jmath + v$ is 10.
\end{proposition}

Now we apply Proposition 
\ref{prop:zero order mod 3 of theta}
and Proposition \ref{prop:degree is 10} to compute the multiplicities of
zeros of 
$\Theta (v, p)=\Theta(\frac{1}{2}\newfont I +v+\int_{\sigma (p)}^p\phi)$.
For an element $v=\xi_1v(p_1) + \cdots + \xi_5 v(p_5)$
the Hamming distance $dis(v)$ of $v$ is defined by
$$
dis(v) = \# \{i \mid \xi_i \neq 0 \}.
$$
Note that $q(v) = dis(v)$ ( mod 3).
If $dis(v) = 0,1,5$, then $\Theta (v, p)=
\Theta(\frac{1}{2}\newfont I +
v+\int_{\sigma (p)}^p\phi)$ is identically zero.
In fact, for example, if $dis(v)=0$, 
then $\Theta (v,p)=0$ for $p = p_1, \dots ,p_5, \sigma (p_1), 
\dots , \sigma (p_5)$. 
Moreover if $p= p_0, p_\infty$, $\Theta (v,p)=0$.
Since the degree of the restriction of theta divisor is 10, 
$\Theta (v,p)$ is identically zero.
For $v$ such that  $dis(v)=2,3,4$, we have the following table.

The case $dis(v)=2$, for example $v=v(p_1)+v(p_2)$,
$$
\matrix
\text{point}	
 & p_1 & p_2 & p_3 & p_4 & p_5 & p_0 \\
\text{order}
& 2 & 2 & 0 & 0 & 0 & 2 \\
\\
\text{point}	
 & \sigma(p_1) & \sigma(p_2) & \sigma(p_3) & 
\sigma(p_4) & \sigma(p_5) & p_\infty \\
\text{order}
& 1 & 1 & 0 & 0 & 0 & 2 .\\
\endmatrix
$$

The case $dis(v)=3$, for example $v=v(p_1)+v(p_2)+v(p_3)$,
$$
\matrix
\text{point}	
 & p_1 & p_2 & p_3 & p_4 & p_5 & p_0 \\
\text{order}
& 0 & 0 & 0 & 1 & 1 & 0 \\
\\
\text{point}	
 & \sigma(p_1) & \sigma(p_2) & \sigma(p_3) & 
\sigma(p_4) & \sigma(p_5) & p_\infty \\
\text{order}
& 2 & 2 & 2 & 1 & 1 & 0 . \\
\endmatrix
$$

The case $dis(v)=4$, for example $v=v(p_1)+v(p_2)+v(p_3)+v(p_4)$,
$$
\matrix
\text{point}	
 & p_1 & p_2 & p_3 & p_4 & p_5 & p_0 \\
\text{order}
& 1 & 1 & 1 & 1 & 2 & 1 \\
\\
\text{point}	
 & \sigma(p_1) & \sigma(p_2) & \sigma(p_3) & 
\sigma(p_4) & \sigma(p_5) & p_\infty \\
\text{order}
& 0 & 0 & 0 & 0 & 2 & 1 .\\
\endmatrix
$$
Since the degree of the function $\Theta (v, p)$ is 10,
the function $\Theta (v,p)$ never vanishes
outside of $\Sigma$ for $v$ such that $dis(v) = 2,3,4$.

\subsection{Rational map $\Cal M_{ab} \to \bold P^{79}$ defined by the theta
constants}
\label{subsec:rational map by theta}
We define a rational map $\Cal M_{ab}\to \bold P^{79}$ 
by using theta constants.
Let us define a subset $S$ of $\bold F_3^5$ by
$$
S=\{ v \in \bold F_3^5-\{0 \} \mid q(v)=0 \}.
$$
The subset of $S$ consisting of $v$ such that
$v\cdot r=0$ (resp. $v\cdot r  \neq 0$) is denoted by $S_r$
(resp. $S_{\bar r}$), where $r=(1,1,1,1,1)$.
Note that $\# S = 80$, $\# S_r =20$ and $\# S_{\bar r}=60$.
We define the theta characteristic
$\Theta_v$ indexed by $v =(v^{(1)}, \dots ,v^{(5)})\in S$ as follows.
Let $\tilde v$ be an element of $\bold Z^5$ such that $\tilde v \equiv v$
(mod $3$). Let 
\begin{align*}
& r_1=(-1,-1,-1,-1,-1),\  r_2=(-1,1,1,1,1), \ r_3=(0,-1,1,0,0) \\
& r_4=(0,0,0,-1,1),\  r_5=(0,1,1,-1,-1), 
\end{align*}
and define $\tilde \beta = (r_1\cdot \tilde v, \dots ,r_5 \cdot \tilde v)$.
Then we have 
$$
(r_1\cdot \tilde v)\bar\delta_1 + \cdots + (r_5\cdot \tilde v)\bar\delta_5 =
v^{(1)}v_1+ \cdots + v^{(5)}v_5 \in \frac{1}{1-\rho}L/L,
$$
where $\bar\delta_i = \frac{1}{1-\rho^2}B_i$.
\begin{lemma}
\label{lem:indep of lift}
Set 
$$
\tilde m = \frac{1}{2}(\bold 1, \bold 1) + 
\frac{1}{3}(-\tilde \beta H, \tilde \beta),
$$
then the theta constant $\Theta_{\tilde m}^3$ is independent of
the choice of a lifting $\tilde v$ of $v$ to $\bold Z^5$.
\end{lemma}
\begin{proof}
If we replace $\tilde v$ by $\tilde v + 3 h$ with $h \in \bold Z^5$,
then $\tilde\beta$ and $\tilde m$ are replaced by $\tilde\beta+ 3 g$ 
and $\tilde m + (-gH, g)$ respectively,
where $g=(r_1\cdot h, \dots ,r_5\cdot h)$.
By applying the quasi-periodicity of theta functions:
$$
\Theta_{(m',m'')+(p,q)} = \bold e (m'\cdot {^tq})\Theta_{(m' ,m'')}
$$
for theta constants to
$m'=\frac{1}{2}\bold 1-\frac{1}{3}\beta H$,
$m''=\frac{1}{2}\bold 1+\frac{1}{3}\beta$,
$p= -gH$ and $q=g$, we have
\begin{align*}
\Theta_{(m',m'')+(p,q)}^3 & =
\bold e(3 \cdot m'\cdot g)\Theta_{(m',m'')}^3 \\ 
& = 
\bold e(3 \cdot (\frac{1}{2}(\bold 1, \bold 1)-\frac{1}{3}\beta H)\cdot ^tg)
\Theta_{(m',m'')}^3 \\
& =
\bold e(\frac{3}{2}(\sum_{i=1}^5r_i \cdot h) - \beta H\cdot ^tg)
\Theta_{(m',m'')}^3. 
\end{align*}
Since the entries of $\sum_{i=1}^5 r_i$ are even, 
$\Theta_{\tilde m}^3$ is independent of the choice of
$\tilde v$.
\end{proof}
The theta constant $\Theta_{\tilde m}^3$ is denoted by $\Theta_v^3$.
Note that $\Theta_v$ is well defined only up to the multiplication
of $\mu_3$.
The isomorphism $\oplus_{i=1}^5 \bar \delta_i \simeq \oplus_{i=1}^5 v_i$
defined by the equality (\ref{eq:base bi and vi}) induces a morphism
$\oplus_{i=1}^5\frac{1}{1-\rho}
\bold Z[\rho ]B_i \to \oplus_{i=1}^5\bold F_3v_i$ and
a homomorphism $\pi : U(4,1) \to O(5, \bold F_3)$.
We have the following proposition.
\begin{proposition}
\label{prop:WE6 action on theta const}
\begin{enumerate}
\item
For $v \in S$, we have $\Theta_{-v}^3= -\Theta_v^3$.
\item
The quotient $\Theta_v^3(\tau)/ \Theta_w^3(\tau)$
comes to be a rational function on 
$\Cal M_{ab} = \Gamma (1-\rho)\backslash\bold B_4$.
\item
For $v \in S$ and $g \in O(4,1)\subset U(4,1)$ such that 
$(\bold 1\cdot g -\bold 1){^t\bold 1}\equiv 0$ (mod $4$),
we have
$\Theta_{v}^3(\iota(g)(\tau))= c(g) \Theta_{v\pi(g)}^3(\tau)$.
Here $c(g)$ is a constant depending only on $g\in O(4,1)$.
As a consequence we have
\begin{equation}
\label{eqn:equivariance}
\frac{\Theta_v^3}{\Theta_w^3}(\iota(g)(\tau))=
\frac{\Theta_{v\pi(g)}^3}{\Theta_{w\pi(g)}^3}(\tau)
\end{equation}
for $v, w \in S$.
\end{enumerate}
\end{proposition}
\begin{remark}
For any $v \in \bold F_3^5$ such that $q(v)=2$, 
we can choose a lifting $\tilde v$
in $\frac{1}{1-\rho^2}\oplus_{i=1}^5 B_i \bold Z$ such that
$h((1-\rho^2)\tilde v)=-2$. For example, we can choose
$\frac{1}{1-\rho^2}(B_1+B_2-2B_5)$ and 
$\frac{1}{1-\rho^2}(4B_1+3B_2+3B_3+6B_5)$ as liftings of 
$\bar\delta_1+\bar\delta_2+\bar\delta_5$
and $\bar\delta_1$, respectively. The reflection for the root
$(1-\rho^2)\tilde v$ satisfies the condition of 3 of Proposition 
\ref{prop:WE6 action on theta const}.
Therefore for any $g \in PO(5, \bold F_3)$, 
there exists an element $\tilde g \in O(4,1)$
such that the equation (\ref{eqn:equivariance}) holds.
\end{remark}
\begin{proof}
1. We choose $-\tilde v$ as a lifting of $-v$. 
For a vector $m \in \bold Q^{10}$, we have $\Theta_{-m}=\Theta_{m}$
for theta constants. Therefore we have
\begin{align*}
\Theta_{-v}  & = 
\Theta_{\frac{1}{2}(\bold 1, \bold 1) +\frac{1}{3}(\tilde\beta H, -\tilde\beta)} \\
& =
\Theta_{- \frac{1}{2}(\bold 1, \bold 1) +\frac{1}{3}(-\tilde\beta H, \tilde\beta)} \\
&=
\Theta_{\frac{1}{2}(\bold 1, \bold 1) +
\frac{1}{3}(-\tilde\beta H, \tilde\beta)+(\bold 1, \bold 1)}. 
\end{align*}
We use the quasi-periodicity for theta constants.
The equality
$$
\bold e (3 \cdot (\frac{1}{2}\bold 1 -\frac{1}{3}\beta H)\cdot ^t\bold 1)
=-1,
$$
yields $\Theta_{-v}^3=-\Theta_{v}^3$.

2. 
It is enough to prove that $\Theta_v(\tau)^3/\Theta_w(\tau)^3$ 
is invariant under the action of generators of $\Gamma (1-\rho)$.
The group $\Gamma (1-\rho)$ is generated by complex reflections, 
(c.f.\cite{ACT})
and the invariance under reflections can be proved similarly as in \cite{Ma}, 
Proposition 5.4,
\cite{Shi}, Lemma 4.1.

3. We apply the transformation formula to compute
$\Theta_m (\iota (g)\tau)$,
where $m={\frac{1}{2}(\bold 1, \bold 1)+
\frac{1}{3}(-\tilde\beta H ,\tilde\beta )}$.
We use the notations $m^{\#}, z^{\#}$ of Proposition 
\ref{prop:theta functional equation}
By the definition of the embedding $O(4,1) \to Sp(10,\bold Z)$
given in (\ref{eq:group embedding}), we have
$$
m^\# = (\frac{1}{2}(\bold 1, \bold 1)+\frac{1}{3}(-\tilde\beta H ,\tilde\beta ))\iota (g)
= \frac{1}{2}(\bold 1 H g H, \bold 1 g) 
+\frac{1}{3}(-\tilde\beta g H,\tilde\beta g).
$$
Since we can choose $\tilde v \cdot g$ as a lifting of $v \cdot \pi (g)$,
$\Theta_{v\pi (g)}$ is equal to 
$\Theta_{m_1}$, where 
$$
m_1=\frac{1}{2}(\bold 1, \bold 1) +\frac{1}{3}(-\tilde\beta g H,\tilde\beta g).
$$
To use the formula (\ref{eq:theta trans}),
we write $m_1=(m', m'')$, $m^\#=m_1+(p,q)$.
Then we have
$$
m'\cdot\;^tq =(\frac{1}{2}\bold 1 +\frac{-1}{3}\tilde\beta g H)
\cdot\;^t (\frac{1}{2}\bold 1 g -\frac{1}{2}\bold 1 ).
$$
Since $\bold 1 g \equiv \bold 1$ ( mod $2$), $3\cdot m' {^tq}$ is an
integer if $(\bold 1 g -\bold 1){^t\bold 1}\equiv 0$ ( mod $4$).
\end{proof}
By the second part of Proposition \ref{prop:WE6 action on theta const},
$\Theta = (\Theta_v)$ defines a rational map
\begin{equation}
\label{eqn:def of theta}
\Theta = (\Theta_v)_{v \in S} : \Cal M_{ab} \to \bold P^{79}.
\end{equation}

\section{80 polynomials and theta constants}
\label{sec:cross ratio}
\subsection{Moduli space of cubic surfaces and the action of $W(E_6)$.}
\label{subsec:three moduli}
We introduce two moduli spaces: the moduli space $\Cal M_{cs}$ of
cubic surface and the moduli space $\Cal M_{6pts}$ of ordered six points on the 
projective plane.

\begin{definition}
\label{def:moduli cs}
The set $\{ (X, \Psi_{cs})\}$ of marked cubic surface is denoted by $\Cal M_{cs}$.
It is equipped with a structure of algebraic variety.
\end{definition}
\begin{remark}
An element $g$ of the automorphism group $Aut (\Gamma_{std})$
of the standard graph $\Gamma$
acts on $\Cal M_{cs}$ from the left by 
$(X,\Psi_{cs}) \mapsto (X, g\circ \Psi_{cs})$.
\end{remark}
\begin{definition}
\label{def:moduli 6pts}
We define an open subset $\Cal M_{6pts}$ of the moduli space of
ordered 6 points in $\bold P^2$ as
$$
\Cal M_{6pts} =\{ p=(P_1,\dots ,P_6) \in ((\bold P^2)^6/Aut(\bold P^2)) \mid
\begin{minipage}{3.5cm}
any 3 points are not \\
\noindent
collinear and 6 points \\
\noindent
are not coconic.
\end{minipage}
\},
$$
where the group $Aut(\bold P^2)$ acts diagonally on $(\bold P^2)^6$.
\end{definition}
For an element $(X, \Psi_{cs})$, 
the contraction $Z$ of $X$ by the ($-1$)-lines
corresponding to $e_1, \dots ,e_6$ is known to be isomorphic to
$\bold P^2$. Therefore the images $P_1, \dots , P_6$ of $e_1, \dots ,e_6$
in $Z$ gives an element
of $\Cal M_{6pts}$.  
This correspondence defines a morphism form $\Cal M_{cs}$
to $\Cal M_{6pts}$, which is known to be isomorphic.

The group $Aut(\Gamma_{std})$
acts on $\Cal M_{6pts}$ via this isomorphism. 
Here we describe the action of $Aut(\Gamma_{std})$ on $\Cal M_{6pts}$ for the
later use.
The involution $r_{ij}$
maps $(P_1, \dots ,P_i, \cdots ,P_j, \dots ,P_6) \mapsto
(P_1, \dots ,P_j, \cdots ,P_i, \dots ,P_6)$. 
Let $b_1 :\tilde P \to \bold P^2$ be the blowing up of $\bold P^2$
with the center $P_1\cup P_2\cup P_3$ and $b_2:\tilde P \to \bold P^2$
be the contraction of the strict transforms of $L_{23}, L_{31}$ and $L_{12}$
to points $P_1', P_2'$ and $P_3'$.
The composite rational map $b_2 \circ b_1^{-1}$ is denoted by $Q_{123}$.
The involution $r_{123}$ 
corresponding to the map
$$
(P_1, \dots ,P_6) \mapsto (P_1', P'_2, P'_3, Q_{123}(P_4),
 Q_{123}(P_5), Q_{123}(P_6))
$$
can be described as follows.
If four points $P_1, \dots ,P_4$ are generic position, 
using the action of $GL(3, \bold C)$, we can
normalize them as $(1:0:0), (0:1:0), (0:0:1)$, and $(1:1:1)$.
If we write the points $P_5, P_6$ as
$(1:x_1:x_3), (1:x_2: x_4)$, then the image of this
$(P_1, \dots ,P_6)$ under $r_{123}$ is 
$$
((1:0:0),(0:1:0),(0:0:1),(1,1,1),(1,x_1^{-1}, x_3^{-1}),
(1,x_2^{-1}, x_4^{-1})).
$$

We define a morphism $per$ from $\Cal M_{cs}$ to $\Cal M_{ab}$ after 
Allcock-Carlson-Toledo.
Let $(X, \Psi_{cs} : \Gamma (X) \to \Gamma_{std})$ be a marked
cubic surface.  The line $\Psi_{cs}^{-1}(e_i)$, $\Psi_{cs}^{-1}(c_i)$ and
$\Psi_{cs}^{-1}(l_{ij})$ are denoted as $E_i$, $C_i$ and $L_{ij}$
respectively. Let $Y$ be the $\mu_3$-covering of $\bold P^3$
branching along $X$.  Then $v_i = \Lambda ([E_i]-[L_{i6}])$
($i= 1, \dots ,5$) forms a orthonormal basis of $J(Y)_{1-\rho}$
and this basis defines a level $(1-\rho)$-structure
$\Psi_{ab} : \bold F_3^5 \to J(Y)$. With the natural polarization on $J(Y)$,
$(J(Y), \iota , \Psi_{ab})$ is an element of $\Cal M_{ab}$.
\begin{definition}
\label{def:period map}
The above correspondence 
$$(X, \Psi_{cs} )
\mapsto (J(Y), \iota, \Psi_{ab})
$$
defines a morphism from $\Cal M_{cs}$ to $\Cal M_{ab}$.
This morphism is called the period map for cubic surfaces and
denoted by $per$.
\end{definition}

\begin{proposition}
\label{prop:compat moduli action}
The actions of $Aut(\Gamma_{std} )$ and $PO(\bold F_3, 5)$ on $\Cal M_{cs}$ and
$\Cal M_{ab}$ are equivariant via the morphism $per$.
\end{proposition}

\subsection{Projective embedding of $\Cal M_{cs}$ defined by 80 polynomials}
\label{subsec:projective embedding}

Let $\tilde P_i=\;^t(P_{i1},P_{i2},P_{i3})$ be an element in $\bold C^3$.
For each element $v \in S$, we attach relatively invariant
polynomials $Z_v$ on $(\bold C^3)^6$ under the action of 
$(g,t)\in GL(3,\bold C) \times (\bold C^\times )^6$
defined by
$\tilde p =(\tilde P_1, \dots ,\tilde P_6) \mapsto 
(t_1 g \tilde P_1,\cdots ,t_6 g \tilde P_6)$,
where $t=(t_1, \dots, t_6)$.

The determinant $\det (\tilde P_i,\tilde P_j,\tilde P_k)$ 
is denoted by $D_{ijk}$.
We define 
$$
\tilde P_i^{(2)}=\;^t(P_{i1}^2,P_{i2}^2,P_{i3}^2,
P_{i1}P_{i2},P_{i2}P_{i3},P_{i3}P_{i1})
$$ 
and
$D_{1\cdots 6}=\det (\tilde P_1^{(2)}, \dots ,\tilde P_6^{(2)})$.
Let $i,j,k,l,m$ be distinct elements in $\{ 1,2,3,4,5 \}$.
For $v_i+ v_j + v_k \in S$, we define a polynomial $Z_{v_i+v_j+v_k}$ by
$$
(-1)^{i+j+k}Z_{v_i+v_j+v_k} = D_{ijk}D_{lm6}D_{1\cdots 6},
$$
if $i < j <k$, $l< m$. In the same way,
we define a polynomial $Z_{v_k+v_l-v_m}$ by
$$
Z_{v_k+v_l-v_m} = 
D_{ikl}D_{jkl}D_{km6}D_{lm6}D_{mij}D_{6ij}. 
$$
For general $v \in S$, we define $Z_v$ by the rule $Z_{-v}= -Z_v$.
If $\tilde P_i$ is not zero, $\tilde P_i$ determines a point in $\bold P^2$
and it is denoted by $P_i$. Note that $Z_v$'s are not zero
if $p = \{P_1,\dots ,P_6\} \in (\bold P^2)^6$ 
is a point of $\Cal M_{6pts}$.  Therefore
$(Z_v)_{v\in S}$ defines a morphism to the projective space $\bold P^{79}$:
\begin{equation}
\label{eqn:coble embedding}
Z=(Z_v)_{ v \in S} :\Cal M_{6pts} \to \bold P^{79}.
\end{equation}
According to \cite{C}, \cite{DO} and \cite{Y}, this projective morphism
is actually an embedding. 
Via the isomorphism from $\Cal M_{cs}$ to $\Cal M_{6pts},$ 
we regard the morphism $Z$ as that from $\Cal M_{cs}$ to $\bold P^{79}$. 
The subgroup in $O(\bold F_3,5)$ generated by $R_{ij}$ ($i < j$)
is isomorphic to $\frak S_6$, and $R_{ij}$ corresponds
to the transposition $(ij)$ via this isomorphism.

and it is identified with $\frak S_6$.
An element $g$ in $\frak S_6$ acts on $(\bold C^3)^6$
by 
$$
\tilde P =(\tilde P_1, \dots ,\tilde P_6) \mapsto 
g \cdot \tilde P = (\tilde P_{\sigma^{-1}(1)}, \dots ,\tilde P_{\sigma^{-1}(6)}).
$$
Note that this action is compatible with the action of $Aut (\Gamma_{std})$
on $\Cal M_{6pts}$.
\begin{proposition}
\label{prop:action of S6 on poly}
For an element $g \in \frak S_6$, we have
$$
Z_{v}(g\cdot \tilde P) = Z_{vg}(\tilde P).
$$
\end{proposition}
\begin{proof}
It is enough to prove the proposition for
$R_{12},\dots ,R_{56} \in \frak S_6$. 
Note that if $v \in S_r$ (resp. $v \in S_{\bar r}$) then
$v \cdot g \in S_r$ (resp. $v \in S_{\bar r}$).
Therefore we can check that the proposition by the definition of $Z_v$.
\end{proof}
For the statement of the main theorem, we recall notations.
\begin{itemize}
\item
$\Cal M_{cs}$: The moduli space of the marked cubic surfaces. See Definition
\ref{def:moduli cs}.
\item
$\Cal M_{ab}$: The moduli space of 5-dimensional abelian varieties with 
actions of $\mu_3$ of type $4\chi \oplus \bar\chi$ and level 
$(1-\rho)$-structures.
See Definition \ref{def:moduli ab}.
\item
$Z$: The morphism from $\Cal M_{cs}$ to $\bold P^{79}$ defined by 80 polynomial
after Coble.
See (\ref{eqn:coble embedding}) in \S \ref{subsec:projective embedding}.
\item
$per$: The period map for the cubic surfaces after Allcock-Carlson-Toledo.
See Definition \ref{def:period map}.
\item
$\Theta$: The morphism form $\Cal M_{ab}$ to $\bold P^{79}$ defined by 80 theta
constants. See (\ref{eqn:def of theta}) in \S 
\ref{subsec:rational map by theta}.
\end{itemize}
The main theorem of this paper is the following. 
\begin{theorem}[Main Theorem]
\label{main theorem}
The following diagram commutes.

\setlength{\unitlength}{0.75mm}
\begin{picture}(160,45)(-40,-5)
\put(0,0){$\Cal M_{ab}$}
\put(0,30){$\Cal M_{cs}$}
\put(60,15){$\bold P^{79}$}
\put(6,15){$per$}
\put(30,27){$Z$}
\put(30,3){$\Theta$}
\put(3,25){\vector(0,-1){20}}
\put(16,4){\vector(4,1){40}}
\put(16,28){\vector(4,-1){40}}
\end{picture}

In other words, $(Z_v)_{v\in S}=(\Theta_v)_{v\in S}$ in $\bold P^{79}$.
In particular, $per$ is an isomorphism from $\Cal M_{cs}$ to
the open set $U$ of $\Cal M_{ab}$ defined by $\cap_{v\in S}\{ \Theta_v \neq 0\}$.
\end{theorem}
Note that the second statement is announced in \cite{ACT}.

\subsection{Proof of Theorem \ref{main theorem}}
\label{subsec:proof of the main theorem}
\begin{proposition}
\label{proportionality}
Via the morphism $per :\Cal M_{6pts}\to \Cal M_{ab}$, we have 
\begin{equation}
\label{eq:ratio}
\frac{\Theta_{v_1+v_2+v_3}^3}{\Theta_{v_1+v_2-v_3}^3}
=c\cdot
\frac{Z_{v_1+v_2+v_3}}{Z_{v_1+v_2-v_3}},
\end{equation}
where $c$ is a 6-th root of unity.
\end{proposition}
\begin{proof}
In this proof, for two non zero functions $f,g$
we write $f \approx g$ if there exists a 6-th root of unity $c$
such that $f= c\cdot g$.
To prove the equality (\ref{eq:ratio}), it enough to show this
on the open set of $\Cal \Cal M_{cs}$ corresponding to cubic surfaces with
no Eckardt points.
We use the normal form of 6 points as in \S\ref{subsec:normal form};
$$
P_i=(1:a_i^2:a_i) \text{  for $i=1, \dots ,5$},\quad
P_6=(0:0:1).
$$
Since $D_{ijk}\approx (a_i-a_j)(a_j-a_k)(a_k-a_i)$
for $i,j,k \neq 6$, $D_{ij6}\approx a_i^2- a_j^2$ and
$D_{1\cdots 6} \approx \prod_{1\leq i < j \leq 5}(a_i -a_j)$,
we have
\begin{equation}
\label{eq:quot for Y}
\frac{Z_{v_1+v_2+v_3}}{Z_{v_1+v_2-v_3}} \approx
\frac{(a_3-a_2)(a_3-a_1)}
{(a_3+a_2)(a_3+a_1)}.  
\end{equation}
For a point $p=(x,y) \in C$, we define a $\bold C^5$ multi-valued 
function $\iota (p) = \int_{\sigma (p)}^p\phi$. As we mentioned, $\iota (p)$
depends on the path connecting $p_0$ and $p$.
We choose a lifting $\tilde v_i=\frac{1}{3}(-s_i H\tau + s_i)$ of $v_i$
in $\frac{1}{1-\rho}L$. 
Then $\rho (\tilde v_i)=\frac{1}{3}(2 s_i H\tau + s_i)$ and
$\rho^2 (\tilde v_i)=\frac{1}{3}(-s_i H\tau -2 s_i)$.
By the quasi-periodicity for
theta functions, the multivalued meromorphic function
\begin{equation}
\label{eq:a rational function}
f(p) = \frac{\prod_{i=0}^2
\Theta (\frac{1}{2}\newfont I + \rho^i(\tilde v_2+\tilde v_3)+\iota (p))}
{\prod_{i=0}^2
\Theta (\frac{1}{2}\newfont I + \rho^i(\tilde v_2-\tilde v_3)+\iota (p))},
\end{equation}
on $p$ becomes a single valued rational function on $C$.
Using the table of \S\ref{subsec:restricted theta}, we have
the following table on the order of zero of the numerator 
and the denominator at $p_i=(a_i, 0)$ of 
(\ref{eq:a rational function}):
$$
\matrix
\text{point}	
 & p_1 & p_2 & p_3 & p_4 & p_5 & p_0 \\
\text{
order of numerator
}
& 0 & 6 & 6 & 0 & 0 & 6 \\
\text{
order of denominator
}
& 0 & 6 & 3 & 0 & 0 & 6 \\
\\

\text{point}	
 & \sigma( p_1) & \sigma( p_2) & \sigma( p_3) & 
\sigma( p_4) & \sigma( p_5) & p_\infty \\
\text{
order of numerator
}
& 0 & 3 & 3 & 0 & 0 & 6 \\
\text{
order of denominator
}
& 0 & 3 & 6 & 0 & 0 & 6 \\
\endmatrix
$$
As a consequence,
the rational function $f(p)$ is equal to 
$c \cdot \displaystyle\frac{x-a_3}{x+a_3}$, where
$c$ is a constant independent of $p$.

First we evaluate the function $f(p)$ at $p_1$. The value
$f(p_1)=c\cdot \displaystyle\frac{a_1-a_3}{a_1+a_3}$
is equal to
\begin{align}
\label{eq:first limit}
c\cdot \displaystyle\frac{a_1-a_3}{a_1+a_3}=& \frac{\prod_{i=0}^2
\Theta (\frac{1}{2}\newfont I + \rho^i(\tilde v_2+\tilde v_3)+\tilde v_1)}
{\prod_{i=0}^2
\Theta (\frac{1}{2}\newfont I + \rho^i(\tilde v_2-\tilde v_3)+\tilde v_1)} \\
\approx &
\frac{\prod_{i=0}^2
\Theta_{\frac{1}{2}(\bold 1, \bold 1) + 
\rho^i(\tilde v_2+\tilde v_3)+\tilde v_1}}
{\prod_{i=0}^2
\Theta_{\frac{1}{2}(\bold 1, \bold 1) + 
\rho^i(\tilde v_2-\tilde v_3)+\tilde v_1}}\cdot
\bold e (\frac{4}{3}s_2H \tau H ^ts_3).
\notag
\end{align}
Here we used the relation between theta functions and theta constants:
$$
\Theta_{(a,b)} = \Theta (a\tau +b)
\bold e (\frac{1}{2}a \tau^t a + a^t b).
$$

Next we consider the limit $\lim_{x \to \sigma(p_2)}f(p) =
c\cdot \displaystyle\frac{-a_2-a_3}{-a_2+a_3}$.
To compute the limit of theta functions, we choose a path
from $Q_0$ to $p$ such that $\iota (p)$ tends to $-\tilde v_2$
if $p$ tends to $\sigma(p_2)$. Using this path, $h =\tilde v_2+\iota (p)$ 
is a function on
$C$ defined on a neighborhood of $\sigma(p_2)$. 
We can choose a local parameter $u$ of $C$
at $\sigma(a_2)$ such that $u(\sigma(p_2))=0$ and $h(-u) = -h(u)$.
Since the order of zero of 
$$
F(u) = 
\Theta (\frac{1}{2}\newfont I + \rho^i(\tilde v_2+\tilde v_3)+\iota (p))
=\Theta (\frac{1}{2}\newfont I + \rho^i(\tilde v_2+\tilde v_3)-\tilde v_2 +h(u))
$$
at $u=0$ is one, we have $\displaystyle\lim_{u \to 0}\frac{F(u)}{F(-u)}=-1$.
Therefore we have
\begin{align}
\label{eq:limit}
& \lim_{x \to \sigma(p_2)}
\frac{\Theta (\frac{1}{2}\newfont I + \rho^i(\tilde v_2+\tilde v_3)+\iota (p))}
{\Theta (\frac{1}{2}\newfont I + \rho^i(\tilde v_2-\tilde v_3)+\iota (p))} \\
 = &
\lim_{u \to 0}
\frac{\Theta (\frac{1}{2}\newfont I + 
\rho^i(\tilde v_2+\tilde v_3)-\tilde v_2 +h (u))}
{\Theta (\frac{1}{2}\newfont I + \rho^i(\tilde v_2-\tilde v_3)-\tilde v_2 +h (u))} 
\notag
\\
 = &
\lim_{u \to 0}
\frac{\Theta (\frac{1}{2}\newfont I + \rho^i(\tilde v_2+\tilde v_3)-\tilde v_2 +h (u))}
{\Theta (\frac{1}{2}\newfont I + \rho^i(\tilde v_2+\tilde v_3)-\tilde v_2 -h (u)+w)}, 
\notag
\end{align}
where $w=-(\bold 1, \bold 1) + 
2\tilde v_2 -2 \rho^i(\tilde v_2)\in \bold Z^{10}$.
Put $w= (w', w'')$ and
use the equality
$$
\Theta (z +a\tau +b) = \bold e (-\frac{1}{2}a\tau^t a - a^t z)\Theta (z),
$$
for $a, b \in \bold Z^5$, then we have
\begin{align*}
& \Theta (\frac{1}{2}\newfont I + \rho^i(\tilde v_2+\tilde v_3)-\tilde v_2 -h (u)+w) \\
= &
\Theta (\frac{1}{2}\newfont I + \rho^i(\tilde v_2+\tilde v_3)-\tilde v_2 -h (u))
\cdot \bold e(-\frac{1}{2}w'\tau ^t w'-w'{^tz_i}),
\end{align*}
where
$z_i= \frac{1}{2}\newfont I + \rho^i(\tilde v_2+\tilde v_3)-\tilde v_2 -h (u)$.
Therefore the limit (\ref{eq:limit}) is equal to
$$
-\bold e(\frac{1}{2}w'\tau ^t w'+
w'{^t}(\frac{1}{2}\bold 1 + \rho^i(\tilde v_2+\tilde v_3)-\tilde v_2)).
$$
Multiplying the equality (\ref{eq:limit}) for $i=0,1,2$, we have 
\begin{equation}
\label{eq:second limit}
c\cdot \frac{-a_2-a_3}{-a_2+a_3}\approx
\bold e (\frac{4}{3}s_2H \tau H^t s_3).
\end{equation}
By the equality (\ref{eq:first limit}) and (\ref{eq:second limit}), 
$\Theta_{v_1+v_2+v_3}^3/\Theta_{v_1+v_2-v_3}^3$ is equal to
the right hand side of (\ref{eq:quot for Y}).

\end{proof}
The following combinatorial lemma is straight forward.
\begin{lemma}
\label{generated by group actions}
We put $e_0=v_1+v_2-v_3, f_0=v_3+v_5-v_4$ (resp.
$e_0'=v_1+v_2+v_6, f_0'=v_1+v_2+v_3$).
Then for any elements $e, f \in S_{\bar r}$ (resp.
$e', f' \in S_{ r}$), there exist
$g_1, \dots ,g_k \in \frak S_6$
(resp. $g_1', \dots ,g_k' \in \frak S_6$) such that
$e = e_0g_1, f_0g_1= e_0g_2, \dots ,f_0g_k=f$
(resp. $e' = e_0'g_1', f_0'g_1'= e_0'g_2', \dots ,f_0'g_k'=f'$).
\end{lemma}

For elements $v,w \in S$, 
we consider the following equality $EC(v,w,c)$:
\begin{align}
\tag*{$EC(v,w,c)$}
\frac{\Theta_{v}^3}{\Theta_{w}^3} = c\cdot \frac{Z_{v}}{Z_{w}}.
\end{align}
For example by Proposition \ref{proportionality}, 
the equality $EC(v_1+v_2+v_3, v_1+v_2-v_3,c)$
holds.
By Proposition \ref{prop:WE6 action on theta const} and
Proposition \ref{prop:action of S6 on poly},
the equality $E(v,w,c)$ implies $E(vg, wg,c)$
for $g \in \frak S_6$.  
\begin{lemma}
\label{lem:EC for Sr and Sbar}
For elements $v,w \in S_r$, (resp. $v,w \in S_{\bar r}$)
the statement $EC(v,w,1)$ holds.
\end{lemma}
\begin{proof}
By applying $(56)\in \frak S_6$ to the equality $EC(v_1+v_2+v_3, v_1+v_2-v_3,c)$,
we have the equality $EC(v_1+v_2+v_3, v_3+v_5-v_4,c)$.  
Therefore we have
\begin{align*}
\frac{Z_{v_1+v_2-v_3}}{Z_{v_3+v_5-v_4}} &=
c\cdot\frac{Z_{v_1+v_2+v_3}}{Z_{v_3+v_5-v_4}}\cdot
c^{-1}\cdot\frac{Z_{v_1+v_2-v_3}}{Z_{v_1+v_2+v_3}} \\
& =\frac{\Theta_{v_1+v_2+v_3}^3}{\Theta_{v_3+v_5-v_4}^3}\cdot
\frac{\Theta_{v_1+v_2-v_3}^3}{\Theta_{v_1+v_2+v_3}^3} \\
& =
\frac{\Theta_{v_1+v_2-v_3}^3}{\Theta_{v_3+v_5-v_4}^3}. 
\end{align*}
Thus we
have the equality $EC(e_0, f_0,1)$, where $e_0=v_1+v_2-v_3$, $f_0=v_3+v_5-v_4$.
For given $e ,f \in S_{\bar r}$,
we can choose $g_1, \dots ,g_k \in \frak S_6$ such that
$e = e_0g_1, f_0g_1= e_0g_2, \dots ,v_fg_k=f$ by Lemma 
\ref{generated by group actions}.
Using these elements $g_1, \dots , g_k$, we have
$$
\frac{Z_{e}}{Z_{f}}
=
\frac{Z_{e_0g_1}}{Z_{f_0g_1}}
\frac{Z_{e_0g_2}}{Z_{f_0g_2}}
\cdots
\frac{Z_{e_0g_k}}{Z_{f_0g_k}},
$$
and the similar equality for $\displaystyle\frac{\Theta_e^3}{\Theta_f^3}$.
The equality $E(e_0, f_0, 1)$ implies
$E(e_0\cdot g_i,f_0\cdot g_i,1)$ and as a consequence, we have
the equality $E(e, f,1)$ for $e, f \in S_{\bar r}$.

In the same way by applying $(36) \in \frak S_6$ to the equality 
$EC(v_1+v_2+v_3, v_1+v_2-v_3,c)$, we get the equality 
$EC(v_1+v_2+v_6, v_1+v_2-v_3,c)$. By taking quotient, we have
$EC(e_0, f_0,1)$, where
$e_0=v_1+v_2+v_6$, $f_0=v_1+v_2+v_3$. Again using Lemma 
\ref{generated by group actions}
and the same technic as $S_{\bar r}$, we get $EC(e,f,1)$
for $e,f \in S_r$.
\end{proof}
\begin{proof}[Proof of the Main Theorem]  
We can directly check the equality
$$
\frac{Z_{v_3+v_4-v_5}}{Z_{v_1+v_2-v_5}}(r_{123}(x))
=
\frac{Z_{-v_1-v_2-v_4}}{Z_{v_2+v_4-v_5}}(x).
$$
By Proposition \ref{prop:WE6 action on theta const}, we have the equality
$EC(-v_1-v_2-v_4,v_2+v_4-v_5,1)$.
As a consequence we have $EC(v,w,1)$ for all $v, w \in S$ and get
the commutativity of the diagram in the main theorem.

We prove the last part of the theorem.  By \cite{Y}, the map 
$\Cal M_{cs} \to U$ is an isomorphism. Thus the map $per$ is proper 
morphism and the image contains an open dense subset and
it is surjective. By the commutativity of diagram, since $per$ is injective,
it is isomorphism. Thus we get the last statement of the main theorem.
\end{proof}

\subsection{Relations for theta constants}
\label{subsec:relation for theta}
As an application of the main theorem \S\ref{main theorem}, we prove 
identities satisfied by theta constants. The
following linear relation between $Z_v$'s is one of the 
Pl\"ucker relation for
$(3\times 3)$-minors in $(3\times 6)$-matrices :
$$
Z_{v_1}+Z_{v_2}+Z_{v_3}+Z_{v_4}=0,
$$
where $v_1=(1,1,0,1,0)$, $v_2=-(1,0,0,1,1)$, $v_3=(0,1,1,0,1)$ and
$v_4=-(1,0,1,1,0)$.
Since the elements $\pm v_1, \dots ,\pm v_4$ are characterized as
the set of $S/\{\pm 1\}$ 
vertical to $w_1=(1,1,1,1,1)$ and $w_2=(1,0,0,-1,0)$.
Note that $q(w_1)=q(w_2)=2$, $w_1\cdot w_2 = 0$ 
and the set $\{ v_1, \dots ,v_4\}$ satisfies the condition
$v_i\cdot v_j =1$ if $i \neq j$.
If another representative $\{ v_1', \dots ,v_4'\}$ of 
$\pm v_1, \dots , \pm v_4$ satisfies the same condition, then
$v'_1=v_1, \dots ,v_4' = v_4$ or
$v'_1=-v_1, \dots ,v_4' = -v_4$.
Since on the set
$$
R = \{\text{ unordered pair }(w_1, w_2) 
\in (\bold P(\bold F_3^5))^2 \mid q(w_1)= q(w_2)=2, w_1\cdot w_2=0\},
$$
the group $PO(\bold F_3, 5)$ acts transitively, we have the 
following cubic relations.
\begin{corollary}
\label{cubic relation}
For $\bold w = (w_1, w_2)\in R$, we define $I(\bold w)$ by
$I(\bold w)=\{ v \in S/\{ \pm 1\} \mid v\cdot w_1=v\cdot w_2=0\}$.
Then we have
$$
\sum_{v \in \tilde I(\bold w)}\Theta_v^3=0.
$$
Here we choose a representative 
$\tilde I(\bold w)=\{v_1, \dots , v_4\}$ 
of $I(\bold w)$ such that $v_i\cdot v_j=1$ if $i\neq j$.
\end{corollary}

On the other hand, one can prove the following cubic relation 
for polynomials $\{Z_v\}$ (See \cite{Y}):
$$
Z_{v_1}Z_{v_2}Z_{v_3}=
Z_{w_1}Z_{w_2}Z_{w_3}, 
$$
where
\begin{align*}
& v_1=(1,1,0,1,0),\ v_2=(0,1,1,-1,0),\ v_3=(1,-1,1,0,0) \\
& w_1=-(1,0,0,1,1),\ w_2=(1,0,1,0,-1),\ w_3=(0,0,1,-1,1). 
\end{align*}
If we put $u=(-1,0,1,1,0)$, then the sets 
$V_1=\{  \pm v_1,\pm v_2 ,\pm v_3, \pm u,0\}$ and
$V_2=\{ \pm w_1, \pm w_2,\pm w_3, \pm u,0\}$ are maximal totally isotropic 
subspaces in $\bold F_3^5$. 
To determine the signature of the equality we consider
$\epsilon =\prod_{i,j}(v_i\cdot w_j)$. 
If we change one of the signatures of
$v_1, \dots ,w_3$, then the signature $\epsilon$ changes.
We define the set
%\begin{align*}
$$
Q=\{ \text{ unordered pair } (V_1, V_2)  
\mid  
\text{
\begin{minipage}{6cm}
 $V_1$ and $V_2$ are
maximal totally  
isotropic subspaces of $\bold F_3^5$ 
and $V_1 \cap V_2$ is 
one dimensional $\bold F_3$ subspace.
\end{minipage}
}
\}. 
$$
%\end{align*}
Then the action of $PO(\bold F_3,5)$ on $Q$ is transitive, we have
the following identity of degree 9.
\begin{corollary}
Let $(V_1,V_2)$ be an element of $Q$. 
Choose a representative $S_{V_1}$ and $S_{V_2}$ of 
$\bold P(V_1)-\bold P(V_1\cap V_2)$
and $\bold P(V_2)-\bold P(V_1\cap V_2)$ in $S$.
Then we have the following 
identity
$$
\prod_{v_1\in S_{V_1}}\Theta_{v_1}^3 = \tilde\epsilon (S_{V_1},S_{V_2})
\prod_{v_2\in S_{V_2}}\Theta_{v_2}^3.
$$
Here the lifting of 
$\prod_{v_1 \in S_{V_1}, v_2 \in S_{V_2}}(v_1\cdot v_2)$
to $\{ \pm 1\}$ is denoted by
$\tilde\epsilon (S_{V_1},S_{V_2})$.
 
\end{corollary}
\begin{remark}
Note that the system of equations
\begin{align*}
&\prod_{v_1\in S_{V_1}}Z_{v_1}  = \tilde\epsilon (S_{V_1},S_{V_2})
\prod_{v_2\in S_{V_2}}Z_{v_2} \quad ((V_1,V_2) \in Q)\\
&\sum_{v \in \tilde I(\bold w)}Z_v =0 \quad (\bold w \in R) 
\end{align*}
is a defining system of equations of the closure of $\Cal M_{6pts}$
in $\bold P^{79}$ (see \cite{Y}).
\end{remark}


\begin{thebibliography}{15}
\setlength{\itemsep}{0.05in}%{-.02in}
\bibitem[ATC]{ACT} Allcock, D., Carlson, J.A. and Toledo, D., 
A complex hyperbolic structure for moduli space of cubic surfaces, 
{\sl C.R Acad. Sci. } {\bf 326} (1998), 49--54.
%
\bibitem[AF]{AF} Allcock, D. and Freitag, E., 
Cubic Surfaces and Borcherds Products, preprint (math.AG/0002066).
%
\bibitem[CG]{CG} Clemens, C.H.  and Griffiths, P.A.  
The intermediate Jacobian of the cubic threefold, 
{\sl Ann. Math.  } {\bf 95} (1969), 460--541.
%
\bibitem[C]{C} Coble, A, 
Points sets and allied Cremona transformations I,II and III, 
{\sl Trans. AMS} {\bf 16} (1915), 155--198, {\bf 17} (1916), 345--385 
and {\bf 18} (1917), 331--372. 
%
\bibitem[DM]{DM} Deligne, P. and  Mostow, G. D., 
Monodromy of hypergeometric functions and nonlattice integral monodromy,  
{\sl I.H.E.S. Publ. Math.} {\bf  63} (1986), 5--89. 
%
\bibitem[DO]{DO} Dolgachev, I. and Ortland, D.,   Point sets in projective 
spaces and theta functions, 
{\sl Asterisque.} {\bf 165} (1988).
%
\bibitem[G]{G} van Geemen, B.,  Private note. 
%
\bibitem[H]{H} Hunt, B. The Geometry of some special Arithmetic Quotients, 
{\sl LNM.} {\bf 1637}, Springer, 1996.
%
\bibitem[I]{I} Igusa, J., Theta Functions, Springer, 1972.
%
\bibitem[N]{N} Naruki, I., 
Cross ratio variety  as a moduli space of cubic surfaces,
{\sl Proc. London Math. Soc.} {\bf 45 no. 3} (1982), 1--30.
%
\bibitem[Ma]{Ma} Matsumoto, K., 
Theta constants associated with the triple coverings of 
the complex projective line branching at six points, 
preprint.
%
\bibitem[Mo]{Mo} Mostow, G. D., 
Generalized Picard lattices arising from half-integral conditions, 
{\sl I.H.E.S. Publ. Math.} {\bf 63} (1986), 91--106. 
%
\bibitem[Mu]{Mu} Mumford, D, Prym varieties I, 
Contributions to analysis 
(a collection of papers dedicated to Lipman Bers), 
325--350, Academic Press, New York, 1974. 
%
%
\bibitem[Pic]{Pic}
E. Picard, Sur les fonctions de deux variables ind\' ependantes 
analogues aux fonctions modulaires,  
Acta Math., {\bf 2} (1883), 114--126.
%
\bibitem[Shi]{Shi} 
Shiga, H., On the representation of Picard modular 
function by $\theta$ constants I-II, 
Publ. RIMS, Kyoto Univ. {\bf 24} (1988), 311--360.
%
\bibitem[T]{T} T. Terada, Fonctions hyperg\' eometriques $F_1$ et 
fonctions automorphes I, II, Math. Soc. Japan {\bf 35} (1983), 451--475; 
{\bf 37} (1985), 173--185.
%
\bibitem[Y]{Y} 
Yoshida, M., A $W(E_6)$-equivariant projective embedding of 
the moduli space of cubic surfaces, 
Kyushu University Preprint series in Mathematics 1999-26.
%

\end{thebibliography}
\end{document}